\documentclass[preprint,3p,a4paper]{elsarticle}
\usepackage{mathrsfs}
\usepackage{amsfonts,amsmath,amssymb,amsthm}

\usepackage{latexsym}
\usepackage{tikz}
\usepackage{hyperref}
\usepackage{geometry}

\usepackage{hyperref}
\usepackage{tikz}
\usepackage[all]{xy}
\usepackage[mathscr]{eucal}
\usepackage{soul}
\usepackage{color}

\newtheorem{thm}{Theorem}[section]

\newtheorem{defn}[thm]{Definition}
\newtheorem{prop}[thm]{Proposition}
\newtheorem{lem}[thm]{Lemma}
\newtheorem{coro}[thm]{Corollary}
\newtheorem{cor}[thm]{Corollary}

\newtheorem{rem}[thm]{Remark}
\newtheorem{exmp}[thm]{Example}
\newcommand{\dda}{\mathord{\mbox{\makebox[0pt][l]{\raisebox{-.4ex}{$\downarrow$}}$\downarrow$}}}
\newcommand{\DDa}{\mathord{\mbox{\makebox[0pt][l]{\raisebox{-.4ex}{$\downarrow$}}$\downarrow$}}}

\newcommand{\UUa}{\mathord{\mbox{\makebox[0pt][l]{\raisebox{.4ex}{$\uparrow$}}$\uparrow$}}}
\newcommand{\ovr}{\overrightarrow}
\newcommand{\da}{\mathord{\downarrow}}
\newcommand{\rom}[1]{\rm{\uppercase\expandafter{\romannumeral #1}}}

\newcommand{\ua}{\uparrow \hspace{-2pt}}

  \newcommand{\upp}{\upharpoonleft \hspace{-3pt}}
\begin{document}
\begin{frontmatter}
\title{Free algebras over directed spaces \tnoteref{t1}}
\tnotetext[t1]{Research supported by NSF of China (Nos. 11871353).}

\author{Yuxu Chen}
\ead{yuxuchen@stu.scu.edu.cn}

\author{Hui Kou\corref{cor}}
\ead{kouhui@scu.edu.cn}

\cortext[cor]{Corresponding author}
\address{Department of Mathematics, Sichuan University, Chengdu 610064, China}


\begin{abstract}
Directed spaces are natural topological extensions of dcpos in domain theory and form a cartesian closed category. In order to model nondeterministic semantics, the power structures over directed spaces were defined through the form of free algebras. We show that free algebras over any directed space exist by the Adjoint Functor Theorem. An c-space (resp.\ b-space) can be characterized as a continuous (resp.\ algebraic) directed space. We show that continuous spaces are just all retracts of algebraic spaces by means of topological ideals, which are generalizations of the rounded ideals. Moreover, by categorical methods, we show that the carrier spaces of free algebras over continuous (resp.\ algebraic) spaces are still continuous (resp.\ algebraic) spaces.


  \vskip 2mm

\end{abstract}

\end{frontmatter}

\section{Introduction}
Directed spaces were introduced in \cite{YK2015} as a topological extended framework for dcpos. They were shown to be closely related with these structures in classical domain theory and they form a category that has many fine properties \cite{CK2020,FK2017,XK2020,XK2022,YK2015}. 

\begin{itemize}
\item[(1)]  Posets endowed with the Scott topology, posets endowed with the weak Scott topology, c-spaces  and locally hypercompact spaces \cite{ERNE2009,XK2022} are all directed spaces. 
\item[(2)]  All directed spaces with continuous maps as morphisms form a cartesian closed coreflective full subcategory of the category of topological spaces \cite{YK2015}.
\item[(3)]  Denote $\bf Dcpo$ the category of dcpos and $\bf DTop$ the category of all directed spaces. Define the embedding functor $\mathcal{E} : {\bf Dcpo} \to {\bf DTop}$ as follows: 
$$ \forall D \in {\rm Ob}{(\bf Dcpo)}, \mathcal{E}D = (D,\sigma(D));$$
$$\forall f \in {\rm Mor}{(\bf Dcpo)}, \mathcal{E}f(x) = f(x), \forall x \in D.$$

Then $\mathcal{E}$ preserves products and the exponential objects. In other words, $\bf Dcpo$ can be viewed as an subcategory of $\bf DTop$ and the products and exponential objects in $\bf Dcpo$ agree with those in $\bf DTop $.

\item[(4)]  By defining the $d$-approximation relation on a directed space, c-spaces can be characterized as continuous directed spaces \cite{XK2022}. When $X$ is a dcpo endowed with the Scott topology, the way-below relation on $X$ is equal to the $d$-approximation relation on $X$. A continuous dcpo can be viewed as a special continuous directed space.
\end{itemize} 

Based on these results above, we can  see that directed spaces are very appropriate topological extended framework for dcpos and are potentially to provide denotational semantics for programming languages. In \cite{XK2020}, the notion of directed lower powerspaces were introduced, which is define through the form of free algebras and is analogous to the notion of lower powerdomains in classical domain theory. The concrete structures of directed lower powerspaces were given and naturally the existence of directed lower powerspaces was proved. Recently, the directed upper and convex powerspaces were constructed in a similar way in \cite{XK2020}. Since all the three kind of directed powerspaces were defined  through the form of free algebras, we try to consider general cases. We will define the notion of dtop-algebras and show that the free dtop-algebras over any directed space exist by the Adjoint Functor Theorem.

In \cite{XK2022}, the notion of $d$-approximation on a directed space was introduced to generalized the notion of way-below relation on a dcpo. It was shown that c-spaces and b-spaces can be characterized as continuous and algebraic directed spaces respectively by introducing the notion of $d$-continuity in a similar way with continuity of dcpos.  Denote ${\rm ID}(E)$ the set of all ideals of $E$. It is well known that a continuous dcpo $E$ can be characterized by that the supremum map $\bigvee : {\rm ID}(E) \to E$ has a left adjoint left inverse \cite{CON03,KO95}, the so-called way-below map $\DDa : E \to ID(E)$. The way-below ideal $\DDa x$ of $x$ is the smallest directed lower segment with supremum $x$. By constructing a canonical algebraic space ${\rm I}_T(X)$ for
each $T_0$ space $X$, we will show that $X$ is continuous if and only if it is a continuous retract of an algebraic space if and only if the supremum map $\bigvee : {\rm I}_T(X) \to X$ has a continuous lower adjoint $\DDa : X \to I_T(X)$. The classical theorem for continuous domains is a special case of the generalized theorem for continuous spaces.

In \cite{X2022}, is was shown that for a continuous directed space $X$, the directed upper, lower and convex  powerspaces over $X$ are continuous as well. The results can be extended to general free dtop-algebras as well.
By means of the result that a space is continuous iff the supremum map has a continuous lower adjoint $\DDa$, we show that the free algebras over a continuous (resp. algebraic) space are still continuous (resp. algebraic) by a simple categorical method which was originated used in the case of dcpos \cite{KO95}. 

\section{Preliminaries}

We assume some basic knowledge of domain theory, category theory and topology, as in, e.g., \cite{AJ94,CON03}

Let $P$ be a poset. For $A\subseteq P$, we set $\da A=\{x\in P: \exists a\in A,  \ x\leq a\}$,\ $\ua A =\{x\in P: \exists a\in P, \ a\leq x\}$.\ $A$ is called a lower or upper set, if $A= \da A$ or $A= \ua A$ respectively. For an element $a\in P$, we use $\da a$ or $\ua a$ instead of $\da\{a\}$ or $\ua \{a\}$, respectively. 

For $x,y\in P$, we say that $x$ is way-below $y$, denoted by $x\ll_s y$, if for any directed subset $D$ of $P$ with the supremum $\bigvee D$ (or $\sup D$) existing, $y\leq \bigvee D$ implies $x\leq d$ for some $d\in D$.\ We call $P$ a continuous poset if $\{a\in P: a\ll_s x\}$ is directed and has $x$ as its supremum for all $x\in P$.

Topological spaces will always be supposed to be $T_0$. For a topological space $X$, its topology is denoted by $\mathcal{O}(X)$ or $\tau$. The partial order $\sqsubseteq$ defined on $X$ by $x\sqsubseteq y \Leftrightarrow  x\in \overline{\{y\}}$
is called the  specialization order, where $\overline{\{y\}}$ is the closure of $\{y\}$.\ From now on, all order-theoretical statements about  $T_0$ spaces, such as upper sets, lower sets, directed sets, and so on,  always refer to  the specialization order $\sqsubseteq $.\ For any net $\xi = (x_i)_{i \in I}$,\ where $I$ is a directed set,\ we write it $(x_i)_I$ for short.\ Given any $x\in X$, $(x_i)_I$ is called converging to $x$,\ denoted by $(x_i)_I  \rightarrow x $,\ if $(x_j)$ is eventually in every open neighborhood of $x$.\

A topological space $X$ is called a c-space if $\forall x \in U \in \mathcal{O}(X), \exists y \in X \text{ s.t. } x \in (\ua y)^\circ \subseteq \ua y \subseteq U$. $X$ is called a b-space if $\forall x \in U \in \mathcal{O}(X), \exists y \in X \text{ s.t. } x \in (\ua y)^\circ = \ua y \subseteq U$.

\vskip 3mm

We now introduce the notion of a directed space.

Let $(X,\mathcal{O}(X))$ be a $T_0$ space. Every directed subset $D\subseteq X$ can be  regarded as a monotone net $(d)_{d\in D}$. Set
$DS(X)=\{D\subseteq X: D \ {\rm is \ directed}\}$ to be the family of all directed subsets of $X$. For an $x\in X$, we denote  $D\rightarrow x$  to mean  that $x$ is a limit of $D$, i.e., $D$ converges to $x$ with respect to $X$. Then the following result is obvious.

\begin{lem}\rm  Let $X$ be a $T_0$ space. For any $(D,x)\in DS(X)\times X$, $D\rightarrow x$ if and only if $D\cap U\not=\emptyset$ for any  open neighborhood of $x$.
\end{lem}

Set $DLim(X)=\{(D,x)\in DS(X)\times X:  \   D\rightarrow x \}$ to be the set of all pairs of directed subsets and their limits in  $X$. Let $X$ be a $T_0$ space. A subset $U\subseteq X$ is called directed-open  if for all $(D,x)\in DLim(X)$, $x\in U$ implies $D\cap U\not=\emptyset$. Obviously, every open set of $X$ is directed-open. Set
$\mathcal{D}(X) =\{U\subseteq X: U \ {\rm is} \  {\rm directed}\text{-}{\rm open}\},$
then $\mathcal{O}(X)\subseteq \mathcal{D}(X)$.

\begin{defn}\rm \cite{YK2015} A  topological space $X$ is said to be a directed space if it is $T_0$ and  every directed-open set is open; equivalently, $ \mathcal{D}(X)=\mathcal{O}(X)$.
\end{defn}

In $T_{0}$ topological spaces, the notion of a directed space is equivalent to the monotone determined space defined by Ern\'e \cite{ERNE2009}.\  Given any space $X$, we denote $\mathcal{D}X$ to be the topological space $(X,\mathcal{D}(X))$. The following are some basic properties of directed spaces.

\begin{thm}\rm \cite{YK2015} \label{convergence}
 Let $X$ be a $T_0$ topological space.
\begin{enumerate}
\item[(1)] For all $U\in \mathcal{D}(X)$, $U=\ua U$.
\item[(2)] $X$ equipped with $\mathcal{D}(X)$ is a $T_0$ topological space such that  $\sqsubseteq_d =\sqsubseteq$, where $\sqsubseteq_d$ is the specialization order relative to $\mathcal{D}(X)$.
\item[(3)] For a directed subset $D$ of $X$, $D\rightarrow x$ iff $D\rightarrow_d x$ for all $x\in X$, where $D\rightarrow_d x$ means that $D$ converges to $x$ with respect to the topology $\mathcal{D}(X)$.
\item[(4)] $\mathcal{D}X$ is a directed space.
\end{enumerate}
\end{thm}

\begin{lem}\rm \cite{YK2015}
Let $X,Y$ be two directed spaces. A map $f:X \to Y$ is continuous iff $f$ preserves all limits of directed subsets, i.e., $\forall (D,x) \in DLim(X), (f(D),f(x)) \in DLim(Y)$.

\end{lem}

\begin{lem}\rm \cite{YK2015} \label{dcon}
Let $X$ be a directed space and $Y$ be a $T_0$ space. A map $f:X \to Y$ is continuous iff $f:X \to \mathcal{D}Y$ is continuous.
\end{lem}

Every dcpo endowed with the Scott topology is a directed space.  The notion of a directed space is a very natural topological extensions of dcpos \cite{FK2017,XK2022,YK2015}. Denote the category of directed spaces with continuous maps as morphism by $\bf DTop$. Then $\bf DTop$ is cartesian closed\cite{YK2015}.
From the following statement we see that is a coreflective full subcateogry of $\bf Top_0$, the category of $T_0$ topological spaces. It follows that $\bf DTop$ is also complete and cocomplete.

\begin{thm}\rm  \cite{YK2015} \label{co}
$\bf DTop$ is a coreflective subcategory of $\bf Top_0$.
\end{thm}

\begin{thm}\rm \cite{HER72} The following statements hold.
\begin{itemize}
\item[(1)] Every coreflective subcategory $\bf C$ of $\bf Top$ is cocomplete. 
\item[(2)] Every coreflective subcategory $\bf C$ of $\bf Top$ is complete and the limit in $\bf C$ is the $\bf C$-coreflection for the limit in $\bf Top$.
\end{itemize}

\end{thm}

The above theorem still holds when restricted in  $T_0$ spaces. And in the case of directed spaces, the $\mathbf{C}$-coreflection for a space $X$ is just $\mathcal{D}X$ \cite{YK2015}. Thus, we gain the following statement.

\begin{coro}\rm
$\bf DTop$ is a complete category. 
\begin{enumerate}
\item Let $\{X_j\}_J$ be a family of directed spaces.
The categorical products of $\{X_j\}_J$, denoted by $\bigotimes_{j \in J} X_j$, is equal to $\mathcal{D}(\prod_{j \in J} X_j)$ up to isomorphism, where $\prod_{j \in J} X_j$ is the topological product of $\{X_j\}_J$.
\item Let $f,g$ be two continuous maps between directed spaces $X$ and $Y$. Then the equalizer of $f,g$ is $\mathcal{D}E$, where $E = \{x \in X: f(x) = g(x)\}$ is endowed with the subspace topology.
\end{enumerate}

\end{coro}

For any two directed spaces $X,Y$, we also denote $X \otimes Y$ to be the categorical product of $X,Y$ in $\bf DTop$. We use $\ovr{x}$ to denote an element of $\bigotimes_{j \in J} X_j$ and $(\ovr{x})_j$ to denote the $j$-th argument of $\ovr{x}$.

\begin{lem}\rm \label{concon}
Let $\{X_j\}_J$ be a family of directed spaces. A directed subset $D = \{\ovr{x_i}\}_{i \in I}$ of $\bigotimes_{j \in J} X_j$ converges to $\ovr{x}$ in $\bigotimes_{j \in J} X_j$ iff $\{(\ovr{x_i})_j\}_{i \in I}$ converges to $(\ovr{x})_j$ in $X_j$ for every $j \in J$.
\end{lem}
\noindent{\bf Proof.} Assume that  $D = \{\ovr{x_i}\}_{i \in I}$ is a directed subset of $\bigotimes_{j \in J} X_j$, and $\ovr{x} \in \bigotimes_{j \in J} X_j$. Since $\bigotimes_{j \in J} X_j = \mathcal{D}(\prod_{j \in J} X_j)$, by Theorem \ref{convergence}, $D \to x$ relative to $\bigotimes_{j \in J} X_j$ iff $D \to x$ relative to  $\prod_{j \in J} X_j$ iff $\{(\ovr{x_i})_j\}_{i \in I}$ converges to $(\ovr{x})_j$ in $X_j$ for every $j \in J$. $\Box$

\begin{lem} \rm \label{prod way}
Let $X,Y$ be two directed spaces. Then $(x_1,y_1) \ll (x_2,y_2)$ in $X \otimes Y$ iff $x_1 \ll x_2$ in $X$ and $y_1 \ll y_2$ in $Y$.
\end{lem}
\noindent{\bf Proof.} Assume that $(x_1,y_1) \ll (x_2,y_2)$ in $X \otimes Y$. Given any directed subset $D \subseteq X$ with $D \to x_2$, consider directed subset $\{(d,y_2)\}_{d \in D} \subseteq X \otimes Y$. Then, $\{(d,y_2)\}_{d \in D} \to (x_2,y_2)$ by Lemma \ref{concon}. There exists some $d_0 \in D$ such that $(x_1,y_1) \leq (d_0,y_2)$, then $x_1 \leq d_0$. It follows that $x_1 \ll x_2$. The same, we have $y_1 \ll y_2$. Conversely, suppose that $x_1 \ll x_2$ in $X$ and $y_1 
\ll y_2$ in $Y$. Given any directed subset $D = \{(x_i,y_i)\}_{i \in I}$ of $X \otimes Y$ with $D \to (x_2,y_2)$, then $\{x_i\}_{i\in I} \to x_2$ in $X$ and $\{y_i\}_{i \in I} \to y_2$ in $Y$ by Lemma \ref{concon}. There exist $i_1,i_2 \in I$ such that $x_1 \leq x_{i_1}, y_1 \leq y_{i_2}$. Pick $i_0 \geq i_1,i_2$, then $(x_1,y_1) \leq (x_{i_0},y_{i_0})$. Therefore, $(x_1,y_1) \ll (x_2,y_2)$. $\Box$

\begin{lem} \rm 
Let $Y$ be a directed space and $\{X_i\}_{1 \leq k \leq n}$ be a family of finite number of directed spaces and $f$ be a map from $\bigotimes_{1\leq k\leq n} X_k$ to $Y$. Then $f$ is continuous iff it is separately continuous, i.e., continuous at each argument.
\end{lem}
\noindent{\bf Proof.} The necessity is obvious. We only prove the sufficiency. Assume that $f$ is separately continuous. It is easy to see that $f$ is monotone. Given any directed susbet $D = \{\ovr{x_i}\}_{i \in I}$  and element $\ovr{y} = (y_1, y_2, \dots,y_n)$ in $\bigotimes_{1 \leq k\leq n} X_k$ such that $D \to \ovr{y}$, we need only to show that $f(D) \to f(\ovr{y})$ in $Y$. Let $f(\ovr{y}) \in U \in \mathcal{O}(Y)$. Consider the directed subset $D_n = \{(y_1,y_2,\dots,y_{n-1}, (\ovr{x_i})_n)\}_{i \in I}$. By Lemma \ref{concon}, $\{(\ovr{x_i})_n\}_{i \in I} $ converges to $y_n$ in $X_n$. Then $f(D_n) \to f(\ovr{y})$ since $f$ is separately continuous. Thus, there exists some $i_n \in I$ such that $f(y_1,y_2,\dots, (\ovr{x_{i_n}})_n) \in U$. Then consider the directed subset $\{(y_1,y_2,\dots,y_{n-2},(\ovr{x_i})_{n-1} , (\ovr{x_{i_n}})_n)\}_{i \in I}$. The same, there exists some $i_{n-1} \in I$ such that $f(y_1,y_2,\dots,y_{n-2},(\ovr{x_{i_{n-1}}})_{n-1} , (\ovr{x_{i_n}})_n) \in U$. Repeating the process, we finally get that $f((\ovr{x_{i_1}})_1,(\ovr{x_{i_2}})_2,\dots,(\ovr{x_{i_{n-1}}})_{n-1} , (\ovr{x_{i_n}})_n) \in U$. Since $I$ is a directed set, there exists a $i_0 \in I$ such that $i_k \leq i_0$ for all $1 \leq k \leq n$. Then, $f((\ovr{x_{i_1}})_1,(\ovr{x_{i_2}})_2,\dots,(\ovr{x_{i_{n-1}}})_{n-1} , (\ovr{x_{i_n}})_n) \sqsubseteq  f((\ovr{x_{i_0}})_1,(\ovr{x_{i_0}})_2,\dots,(\ovr{x_{i_{0}}})_{n-1} , (\ovr{x_{i_0}})_n) = f(\ovr{x_{i_0}}) \in U$. Thus, $f(D) \cap U \not = \emptyset$. It follows that $f(D) \to f(\ovr{y})$ in $Y$. $\Box$

\vskip 3mm
Next, we introduce the notion of continuous spaces and algebraic spaces.

\begin{defn} \rm 
Let $X$ be a directed space and $x,y\in X$. 
 We say that $x$ $d$-approximates $y$, denoted by $x\ll_d y$, if for any directed subset $D\subseteq X$, $D\rightarrow y$ implies $x\sqsubseteq d$ for some $d\in D$. If $x\ll_d x$, then $x$ is called a $d$-compact element of $X$.
\end{defn}

For any directed space $X$ and $x\in X$, we denote $K_d(X)  = \{x\in X: x\ll_d x\}$ to be the set of all $d$-compact elements of $X$, $\DDa_d x = \{y\in X: y\ll_d x\}$, and $ \UUa_d x = \{y\in X: x\ll_d y\}$. It is easy to see that the $\ll_d$ relation is equal to $\ll_s$ for a dcpo endowed with the Scott topology. The $d$-approximation is a natural extension of way-below relation. Similarly, the notion of a continuous space can be defined and it was found to be equivalent to the notion of a c-space.

\begin{defn} \rm 
A topological space $X$ is called a continuous space if it is a directed space and for any $x\in X$, there exists a directed subset $D \subseteq \DDa_dx$ such that $D \rightarrow x$.
\end{defn}

\begin{lem} \rm \label{topo ideal}
Let $X$ be a $T_{0}$ topological space and $D$ be a directed subset of $X$.\
If $D\subseteq\da x$ and $D\rightarrow x$ for $x\in P$, then $x=\sup D$.

\end{lem}
\noindent{\bf Proof.}
Suppose that $D\subseteq\da x$ and $D\rightarrow x$ for $x\in X$.\ Assume that there exists $y\in X$ such that $D\subseteq \da y$ and $x\not\leq y$.\ Then $x\in P\backslash \da y\in \mathcal{O}(X)$. Thus $D\cap (X\backslash \da y ) \not=\emptyset$, a contradiction.\ Hence,\ $x=\sup D$.\ $\Box$

\vskip 3mm

The following are some basic properties of continuous spaces, which extend these classical results for continuous domains.

\begin{thm} \rm \cite{CK2020,XK2022} \label{d-continuity}
Let $X$ be a continuous space. 
\begin{enumerate}
\item[(1)] For all $x,y\in X$, $x\ll_d y$ implies $x\ll_d z\ll_d y$ for some $z\in X$.
\item[(2)] $\UUa_d x=(\ua x)^{\circ}$ for all $x\in X$. Moreover, $\{\UUa_d x: x\in X\}$ is a base of the topology of $X$.
\item[(3)] For any $x \in X$, $\DDa_d x$ is a
directed set and $\DDa_d x \to x$. Moreover, $x$ is the supremum of $\DDa x$.
\end{enumerate}
\end{thm}

\begin{thm} \rm \cite{XK2022}
A topological space $X$ is a continuous space iff it is a c-space.
\end{thm}

From Theorem \ref{d-continuity}, we see that in a continuous space $x \ll_d y$ iff $y \in (\ua x)^\circ$ and $x$ is a compact element iff $\ua x$ is an open subset. 
Continuous dcpos endowed with the Scott topology can be viewed as special continuous spaces \cite{XK2022,ZHA2021}. From now on, we drop the subscript $d$ in these notations when there is no confusion can arise.
Similar, the notion of an algebraic space can be defined and it is equivalent to the notion of a b-space. 

\begin{defn}\rm
A topological space is called an algebraic space if it is a directed space such that for any $x\in X$, $\da x \cap K(X)$ is a directed subset and converges to $x$.
\end{defn}

Obviously, algebraic spaces are continuous spaces. Finally, we introduce the notion of a basis of a directed space.

\begin{defn}\rm 
Let $X$ be a directed space. A subset $B$ of $X$ is called a basis of $X$ if for any $x \in X$, $\DDa x \cap B$ is a directed subset and converges to $x$.
\end{defn}\rm 

From the definition, we see that a topological space is continuous iff it has a basis.
Let $X$ be an algebraic space. By Lemma \ref{topo ideal}, for a directed subset $D$ of $\da x \cap K(X)$ that converges to $x$, we have $x = \sup D$. Thus, any element of an algebraic space is the supremum and a limit of the set of all compact elements below it, i.e., $K(X)$ is a basis of $X$. Moreover, $K(X)$ is the least basis of $X$.

\begin{prop}\rm
Let $X$ be an algebraic space. $K(X)$ is a basis of $X$. For any basis $B$ of $X$, $K(X) \subseteq B$.
\end{prop}

\noindent{\bf Proof.}
Assume that $B$ is a basis of $X$. For any $x \in K(X)$, we have $\DDa x = \da x$ and $\UUa x = \ua x$. Then $\DDa x \cap B = \da x \cap B$ is a directed subset and converges to $x$. If $x \not\in B$, then $\ua x = \UUa x$ is an open subset containing $x $, which contradicts to that $\da x \cap B $ converges to $x$. Thus, $K(X) \subseteq B$. $\Box$

\begin{thm}\rm 
Let $X$ be a topological space. Then the following conditions are equivalent. 
\begin{enumerate}
\item[(1)] $X$ is an algebraic space.
\item[(2)] $X$ is a directed space and $K(X)$ is a basis of $X$.
\item[(3)] $X$ is a directed space and for any $x \in X$, there exists a directed subset $D \subseteq K(X)$ such that $D 
\to x$ with $\sup D = x$.
\item[(4)] $X$ is a b-space.

\end{enumerate}

\end{thm}
\noindent{\bf Proof.} $(1) \Rightarrow (2) \Rightarrow (3)$ is easily seen from the definition.

$(3)\Rightarrow (4).$ Given any $x \in U \in \mathcal{O}(X) = \mathcal{D}(X)$, supposing that $D \subseteq K(X)$, $D \to x$ and $\sup D = x$, then there exists some $d \in D \subseteq K(X)$ such that $d \in U$ and $x \in \ua d$. We need only to show that $\ua d$ is an open subset of $X$. Let $E$ by any directed subset of $X$ such that $E \to a \in \ua d$, then $E \to d$ as well. Since $d$ is a compact element, there exists some $e \in E$ such that $d \leq e$. Thus, $E \cap \ua d \not = \emptyset$. It follows that $\ua d$ is an directed-open set hence an open set. 

$(4) \Rightarrow (1).$
Suppose that $U \subseteq X$ is a directed-open subset. Given any $x \in U$, let $A_x = \{y \in X: x \in (\ua y)^\circ = \ua y\}$. Given any $d_1,d_2 \in A_x$, then $\ua d_1 \cap \ua d_2$ is an open subset. By definition, there exists a $d_0 \in \ua d_1 \cap \ua d_2$ such that $x \in (\ua d_0)^\circ = \ua d_0$, i.e., $d_1,d_2 \leq d_0 \in A_x$. Thus, $A_x$ is a directed subset of $X$ and $A_x \to x$ by definition of the b-space. Then $A_x \cap U \not = \emptyset$. It follows that there exists some $a \in A_x$ such that $x \in (\ua a)^\circ = \ua a \subseteq U$. Thus, $U$ is an open subset. $X$ is a directed space.

Let $d$ be any element of $X$ such that $(\ua d)^\circ = \ua d$. For any directed subset $D \to d$, there exists some $d^\prime \in D$ such that $d^\prime \in \ua d$, i.e., $d \leq d^\prime \in D$. Thus, $d$ is a compact element of $X$. Conversely, suppose that $d$ is a compact element of $X$ and $d \leq x \in X$. Since $A_x$ is a directed subset and $A_x \to x$, there exists some $a \in A_x$ such that $d \leq a$. Then $x \in (\ua a)^\circ = \ua a \subseteq \ua d$. It follows that $x \in (\ua d)^\circ$ and then $(\ua d)^\circ = \ua d$. Therefore, $\da x \cap K(X) = A_x$ is a directed subset that converges to $x$ for any $x \in X$. $X$ is algebraic.  $\Box$

\vskip 3mm
In the proof of the above theorem, we also in fact proved the following statement.
\begin{prop}\rm \label{basealg}
Let $X$ be an algebraic space. Then, $\ua x = (
\ua x)^\circ $ for each $x \in K(X)$. Moreover, $\{\ua x: x \in K(X)\}$ is a base of the topology of $X$.
\end{prop}

\begin{rem} \rm 
When defining continuous spaces and algebraic spaces, we require them to be directed spaces. Otherwise, it may
not be a c-space or a b-space. The following example illustrates this.

Let $\mathbb{N}$ be the set of natural numbers. Denote $\mathbb{N}^\top$ the flat domain, i.e., the poset with carrier set $\mathbb{N} \cup \{\top\}$ and $x \leq y$ iff $y = \top$ or $x=y$. Consider $\mathbb{N}^\top$ endowed with the upper topology. Then it is easy to verify that for any $x \in \mathbb{N}^\top$, $x \ll_d x$. Thus, $K_d(X) = X$ is a basis of $X$, $\da x \cap K_d(X) \to x$ for any $x \in X$. However, it is not a directed space and hence not a b-space.
\end{rem}

\begin{thm} \rm \cite{XK2022} \label{core compact}
Let $X$ be a directed space. If $X$ is core compact, then for any directed space $Y$, $X \otimes Y = X \times Y$.
\end{thm}

We know that $X$ is a continuous space iff $\mathcal{O}(X)$ is completely distributive \cite{Erne91}. Therefore, continuous spaces are core compact and their finite categorical products in $\bf DTop$ agree with the topological products.

\section{Existence of free algebras over directed spaces}

The notion of directed lower powerspaces \cite{XK2020} was introduced to provide the nondeterministic denotational semantics. It is defined through the from of free algebras with an operation that is idempotent, commutative, associative and deflationary. Recently, the concepts of upper directed  powerspaces and  directed convex powerspaces were defined in a same way \cite{X2022}. Generally speaking, they are free algebras over directed spaces. In this section, we consider the general free algebras over directed spaces and show the existence of free algebras over directed spaces.

\vskip 3mm
Similarly with dcpo-algebras \cite{AJ94}, we define $(\Sigma,\mathcal{E})$-algebras based on directed spaces, called dtop-algebras. A finitary signature $\Sigma = \langle I,\alpha \rangle$ consists of a set $I$ of operation symbols and a map $\alpha : I \to \mathbb{N}$, assigning to each operation symbol a finite arity. We also write $n_i$ for $\alpha(i)$ for any $i \in I$. A $\Sigma$-algebra  $\langle A,\{f_i\}_I \rangle$ is given by a carrier set $A$ and a set $\{f_i\}_I$ of interpretations of the operation symbols, in the sense that for $i \in I$, $f_i$ is a map from $A^{n_i}$ to $A$. For convenience, we also call each $f_i$ an operation.
A homomorphism between two $\Sigma$-algebras $\langle A, \{k_i\}_I \rangle$ and $\langle B ,\{ l_i\}_I\rangle$ is a map $h: A \to B$ which commutes with the operations:
$$\forall i \in I, h(k_i(a_1,\dots,a_{n_i})) = f_i(h(a_1),\dots,h(a_{n_i}))$$
The term algebra over a set $X$ with respect to a signature $\Sigma$ is denoted by $T_\Sigma(X)$. It has the universal property that each map from $X$ to $A$, where $\langle A, \{f_i\}_I \rangle$ is a $\Sigma$-algebra, can be extended uniquely to a homomorphism $\overline{h}: T_\Sigma(X) \to \langle A,\{f_i\}_I \rangle$. Let $V$ be
a fixed countable set of variables.
$T_\Sigma(V)$ are used to encode equations. An inequality $\tau_1 \leq \tau_2$ is said to hold in an algebra
$\langle A,\{f_i\}_I \rangle$ if for each map $h: V \to A$ we have $\overline{h}(\tau_1) \leq \overline{h}(\tau_2)$. The pair $\langle \overline{h}(\tau_1),\overline{h}(\tau_2)\rangle$ 
is also called an instance of the inequality $\tau_1 \leq \tau_2$.

A $ \Sigma$-algebra is called a dtop-algebra if the carrier set is equipped with a topology such that it becomes a directed space, and such that each operation is continuous. A change should be noticed that for an operation $k:X^n \to X$,  $X^n$ means the directed product of $X$ by $n$ times, i.e., $X^n = X \otimes X 
\otimes \dots \otimes X$. We also require the homomorphisms to
be continuous. We let a pair $\langle \tau_1, \tau_2 \rangle \in  \mathcal{E}  \subseteq T_\Sigma(V ) \times  T_\Sigma(V)$ stand for the inequality
$\tau_1 \leq \tau_2$. We denote ${\bf DTop}(\Sigma, \mathcal{E})$ the category of all dtop-algebras with
signature $\Sigma$ that satisfy the inequalities in $\mathcal{E}$, where the morphisms are all the homomorphisms.

\begin{defn}\rm
Let $X$ be a directed space. A dtop-algebra $\langle A,\{f_i\}_I\rangle$ is called the free algebra over $X$ with respect to ${\bf DTop}(\Sigma,\mathcal{E})$ if $A$ is an object of ${\bf DTop}(\Sigma,\mathcal{E})$ and there exists a continuous map $i: X \to A$ such that any continuous map $f: X \to B$, where $\langle B,\{g_i\}_I\rangle$ is an object of ${\bf DTop}(\Sigma,\mathcal{E})$, extends uniquely to a homomorphism $\overline{f}: \langle A,\{f_i\}_I\rangle \to \langle B,\{g_i\}_I\rangle$ such that $\overline{f} \circ i =f$.
$$\xymatrix{
  X \ar[r]^{i} \ar[dr]_{\forall f} &  A \ar[d]^{\overline{f}}   & \langle A,\{f_i\}_I\rangle \ar[d]_{\exists!\overline{f}}  \\
 & B & \langle B,\{g_i\}_I\rangle
                     }$$

\end{defn}

\vskip 2mm

It is well known that for dcpos and topological spaces, the forgetful functor has a left adjoint \cite{AJ94,CON03} and then the free algebras exist. In the following, we will use the Adjoint Functor Theorem \cite{RI} to prove the existence of the free dtop-algebras over directed spaces with respect to any signature and a set $\mathcal{E}$ of inequalities.

\begin{thm}[Adjoint Functor Theorem] \rm  \cite{RI} \label{adjoint}
 Let $U : \mathbf{A \to S}$ be a continuous functor whose domain is locally small and complete. Suppose that $U$ satisfies the following $\bf solution\ set\ condition:$
\begin{itemize}
\item  For every $s \in \mathbf{S}$ there exists a set of morphisms $\Phi_s = \{f_i : s \to Ua_i\}$ so that any $f:s \to Ua$ factors through some $f_i \in \Phi_s $ along a morphism $a_i \to a$ in $\bf A$.

\end{itemize}
Then $U$ admits a left adjoint.
\end{thm}

Let $U$ be the forgetful functor from $\mathbf{DTop}(\Sigma,\mathcal{E})$ to $\bf DTop$. Then the free dtop-algebras over $X$ exist if $U$ has a left adjoint.
Now we check that $\mathbf{DTop}(\Sigma,\mathcal{E})$ satisfies the conditions.

\begin{lem} \rm  \label{left1}
$\mathbf{DTop}(\Sigma,\mathcal{E})$ is a locally small and complete category.
\end{lem}
\noindent{\bf Proof:} The locally small property is obvious. For completeness, we need only to show that it is closed under products and equalizers. 

1. Products. Let $\{\langle X_j,\{k_{ij}\}_I \rangle\}_{j\in J}$ be a family of objects in $\mathbf{DTop}(\Sigma,\mathcal{E})$. Denote $\prod_{j \in J} \langle X_j,\{k_{ij}\}_I\rangle$ as follows:

\begin{enumerate}
\item[(1)] The carrier space of $\prod_{j \in J} \langle X_j,\{k_{ij}\}_I\rangle$ is $\bigotimes_{j \in J} X_j$. 
\item[(2)] For each $i \in I$, the corresponding operation from $(\bigotimes_{j \in J} X_j)^{n_i}$ to $\bigotimes_{j \in J} X_j$, denoted by $k_i$, satisfies 
$$(k_i(\ovr{x_1},\dots,\ovr{x_{n_i}}))_j = k_{ij}((\ovr{x_1})_j,\dots,(\ovr{x_{n_i}})_j)$$
 
\end{enumerate}
That is, $\prod_{j \in J} \langle X_j,\{k_{ij}\}_I\rangle = \langle \bigotimes_{j \in J} X_j, \{f_i\}_I\rangle$. We claim that it is the categorical product of the family $\{\langle X_j,\{k_{ij}\}_I \rangle\}_{j\in J}$ in $\mathbf{DTop}(\Sigma,\mathcal{E})$.

\vskip 2mm
First we show that $\prod_{j \in J} \langle X_j,\{k_{ij}\}_I \rangle$ is a well defined dtop-algebra. Since $\bigotimes_{j \in J} X_j$ is a directed space, we need only to prove that $k_i$ is continuous for each $i \in I$ and inequalities in $\mathcal{E}$ hold. Denote projection $\pi_j :\prod_{j \in J} \langle X_j,\{k_{ij}\}_I\rangle  \longrightarrow (X_j,\{k_{ij}\}_I)$ for each $j \in J$ as follows: 
$$\forall \vec{x} \in \bigotimes_{j \in J} X_j,\ \pi_j(\vec{x}) = (\vec{x})_j.  $$ 
Given any directed subset $D = \{(\ovr{x_{1l}},\dots,\ovr{x_{n_il}}) \}_{l \in L}$ of $(\bigotimes_{j \in J}X_j)^{n_i}$ and $D \to (\ovr{x_1},\dots,\ovr{x_{n_i}}) = a  \in (\bigotimes_{j \in J}X_j)^{n_i}$, we have 
$$\pi_j (k_i(D)) = \{ k_{ij}((\ovr{x_{1l}})_j,\dots,(\ovr{x_{n_{i}l}})_j)\}_{l \in L},$$
and 
$$\{((\ovr{x_{1l}})_j,\dots,(\ovr{x_{n_il}})_j)\}_L \to ((\ovr{x_1})_j,\dots,(\ovr{x_{n_i}})_j) \in X_j^{n_i}$$
by Lemma \ref{concon},
where $\{((\ovr{x_{1l}})_j,\dots,(\ovr{x_{n_il}})_j)\}_L$ is a directed subset of $X_j^{n_i}$ for each $j \in J$.

Since $k_{ij}: X_j^{n_i} \to X_j$ is continuous, then 
$$ k_{ij} (\{((\ovr{x_{1l}})_j,\dots,(\ovr{x_{n_il}})_j)\}_L ) \to k_{ij}((\ovr{x_1})_j,\dots,(\ovr{x_{n_i}})_j )\in X_j.$$
Thus, $ k_i(D) \to k_i(a)$ by definition of $k_i$ and Lemma \ref{concon}, i.e., $k_i :(\bigotimes_{j \in J} X_j)^{n_i} \to \bigotimes_{j \in J} X_j$ is continuous.

\vskip 2mm
Since $ \ovr{x} \leq \ovr{y}$ in $\bigotimes_{j \in J} X_j$ iff $ (\ovr{x})_j \leq  (\ovr{y})_j$ in $X_j$ for each $j \in J$, by the definition of $k_i$, to verify that an inequality holds, we need only to check that it holds in each argument. It follows from that any inequality in $\mathcal{E}$ holds for $(X_j,\{k_{ij}\}_I) $ for any $j \in J$.

\vskip 2mm


Next, we show that $\pi_j$ is a well defined morphism in $\mathbf{DTop}(\Sigma,\mathcal{E})$.  Obviously, $\pi_j: \bigotimes_{j \in J} X_j \longrightarrow X_j$ is continuous. We need only to check the commutativity of $\pi_j$ and all operations $k_i$ for $i \in I$. We have

  $$\pi_j(k_i (\ovr{x_1},\dots,\ovr{x_{n_i}}))  \\
     =\  k_{ij}((\ovr{x_1})_j,\dots,(\ovr{x_{n_i}})_j) \\
     =\  k_{ij}(\pi_j(\ovr{x_1}),\pi_j(\ovr{x_2}),\dots,\pi_j(\ovr{x_{n_i}})). 
    $$
Thus, $\pi_j$ is a homomorphism for each $j \in J$.
\vskip 2mm

Finally, we show $\prod_{j \in J} \langle X_j,\{k_{ij}\}_I\rangle $ is the categorical product of the family of dtop-algebras $\{\langle X_j,\{k_{ij}\}_I \rangle \}_{J}$. Given any $\langle Z,\{ h_i\}_I\rangle$ in  ${\bf DTop}(\Sigma,\mathcal{E})$ and homomorphisms $f_j:\langle Z,\{ h_i\}_I\rangle \to \langle X_j,\{ k_{ij}\}_I\rangle $ for $j \in J$, we define the morphism  $\prod_{j \in J} f_j$ from $\langle Z,\{ h_i\}_I\rangle$ to $\prod_{j \in J} \langle X_j,\{k_{ij}\}_I\rangle $ as follows: 
$$ ((\prod_{j \in J} f_j) (z))_j = f_j(z), \forall z \in Z.$$ Obviously, it is continuous and $ \pi_j \circ (\prod_{j \in J} f_j) = f_j$. We need only to show the commutativity between $ \prod_{j \in J} f_j$ and all the operations $h_i$ for $i\in I$. We have
\begin{align*}
  & ((\prod_{j \in J} f_j) (h_i (z_1,z_2,\dots,z_{n_i})))_j \\
  = &\  f_j(h_i (z_1,z_2,\dots,z_{n_i}))\\
  = &\  k_{ij}(f_j(z_1),f_j(z_2),\dots,f_j(z_{n_i})) \\
  = &\  (k_i((\prod_{j \in J} f_j)(z_1),\dots,(\prod_{j \in J} f_j)(z_{n_i}))))_j. \\
  &
\end{align*}
Thus,  $\prod_{j \in J} \langle X_j,\{k_{ij}\}_I\rangle $ is the categorical product.

\vskip 3mm

2. Equalizers. Let $f,g:\langle X,\{k_i\}_I\rangle \to \langle Y,\{l_i\}_I\rangle$ be two morphisms in ${\bf DTop}(\Sigma,\mathcal{E})$. Define $E = \{x \in X: f(x) =g(x) \}$ be the subspace of $X$. Let  $h_i(a_1,a_2,\dots,a_{n_i}) =  k_i(a_1,a_2,\dots,a_{n_i}) $ for all $i \in I$ and $a_j \in E$ for all $j \leq n_i$. Then $h_i(a_1,a_2,\dots,a_{n_i}) \in E$ since 
\begin{align*}
& f(h_i(a_1,a_2,\dots,a_{n_i})) \\
=\ &  h_i(f(a_1),f(a_2),\dots,f(a_{n_i}))          \\
=\ &  h_i(g(a_1),g(a_2),\dots,g(a_{n_i}))  \\
=\ &  g(h_i(a_1,a_2,\dots,a_{n_i})). \\
&
\end{align*}
We claim that $\langle \mathcal{D}E,\{h_i\}_I\rangle$ is an object in $\mathbf{DTop}(\Sigma,\mathcal{E})$ and it is an equalizer for $f$ and $g$.
\vskip 2mm

First, we show it is well defined. Of course, the inequalities in $\mathcal{E}$ also hold for  $\langle \mathcal{D}E,\{h_i\}_I\rangle$. We need only to show $h_i : (\mathcal{D}E)^{n_i} \to \mathcal{D}E$ is continuous. Given any directed subset $\{(a_{1j},a_{2j},\dots,a_{n_ij})\}_J$ of $(\mathcal{D}E)^{n_i}$ with $\{(a_{1j},a_{2j},\dots,a_{{n_i}j})\}_J \to (a_1,a_2,\dots,a_{n_i}) \in (\mathcal{D}E)^{n_i} $, we have $\{a_{mj}\}_J \to a_m$ relative to $\mathcal{D}E$ for each $1 \leq m \leq n_i$. By Theorem 
\ref{convergence}, the convergence also holds relative to $E$ and then relative to $X$ since $E$ is a subspace of $X$. Thus,  $\{ h_i ( a_{1j},a_{2j},\dots,a_{{n_i}j}) \}_J = \{ k_i ( a_{1j},a_{2j},\dots,a_{{n_i}j}) \}_J \to k_i(a_1,a_2,\dots,a_{n_i}) = h_i(a_1,a_2,\dots,a_{n_i})$ relative to $X$. And naturally the convergence also holds relative to $\mathcal{D}E$. Thus, $h_i$ is continuous, and $\langle \mathcal{D}E,\{h_i\}_I\rangle$ is a well defined object in  $\mathbf{DTop}(\Sigma,\mathcal{E})$.

The embedding map $em : \langle \mathcal{D}E,\{h_i\}_I\rangle \to \langle X,\{f_i\}_I\rangle$ is a morphism, where $\forall a \in E.\ em(a) = a $. The continuity is obvious. The commutativity between $em$ and the operations is easily seen from the definition.

Finally, we show that $\langle \mathcal{D}E, \{h_i\}_I \rangle$ is an equalizer for $f$ and $g$. Obviously, $f \circ em = g \circ em$. Let  $\langle Z,\{r_i\}_I\rangle$ be any object with a morphism $t:\langle Z,\{r_i\}_I\rangle \to \langle X,\{f_i\}_I\rangle $ such that $f \circ t =g \circ t$. Define $ \overline{t} : \langle Z,\{r_i\}_I\rangle \to \langle \mathcal{D}E,\{h_i\}_I\rangle $ as follows: $\forall z\in Z.\ \overline{t} (z) = t(z)$. Since $f \circ t =g \circ t$, we know $t(z) \in E$. $\overline{t}$ is a well defined map. Since $t: Z \to X$ is continuous, $t: Z \to E$ is continuous and then $\overline{t}: Z \to \mathcal{D}E$ is continuous by Lemma \ref{dcon}. The commutativity between $\overline{t}$ and operations is easily seen. Thus, $\overline{t}$ is a morphism. Obviously, $em \circ \overline{t} = t$. The uniqueness of $\overline{t}$ is also obvious. Thus, we get the conclusion that $\langle \mathcal{D}E,\{h_i\}_I\rangle$ with $em : \langle \mathcal{D}E,\{h_i\}_I\rangle \to \langle X,\{f_i\}_I\rangle$ is the equalizer for $f,g$.


\begin{lem} \rm  \label{left2}
The forgetful functor $U : {\bf DTop}(\Sigma,\mathcal{E}) \to {\bf DTop}$ satisfies the solution set condition.
\end{lem} 
\noindent{\bf Proof.}
$U$ is a continuous functor since it preserves products and equalizers. Given any $X$ in ${\bf DTop}$, any $(Y,\{l_i\}_I)$ in ${\bf DTop}(\Sigma,\mathcal{E})$ and any continuous map $i: X \to Y$, let $\overline{X}$ be the ordinal subalgebra, equipped with the induced topology, of $Y$ which is generated by $i(X)$, i.e., the least subset $A$ of $Y$ containing $i(X)$ such that $l_i(A^n) \subseteq A$ for every $i \in I$. (This can be done by iteration. Denote $S_0 = i(X)$ and $S_{n+1} = S_n \cup \{l_i(x_1,\dots,x_{n_i}) : i \in I\  \&\ x_j \in S_n \text{ for } 1 \leq j \leq n_i\}$ for $n \in \mathbb{N}$. Then $\overline{X} = \bigcup_{n \in \mathbb{N}} S_n$.) Then, the cardinal of $\overline{X}$ is bound depending on the cardinal of $X$ and the signature $\Sigma$. We define $(\mathcal{D}\overline{X},\{k_i\}_I)$ the dtop-subalgebra of $(Y,\{l_i\}_I)$ generated by $X$, where $k_i$ is $l_i$ restricted in $\overline{X}$ for each $i \in I$. Then $k_i :(\mathcal{D}\overline{X})^n \to \mathcal{D}\overline{X}$ is continuous by the fact that any directed subset $\{(a_{1j},a_{2j},\dots,a_{nj})\}_J$ converges to $(a_1,a_2,\dots,a_n)$ in $(\mathcal{D}\overline{X})^n$ iff $\{a_{mj}\}_J \to a_m$ in $\overline{X}$ for each $1 \leq m\leq n$,  $l_i$ is continuous and $k_i$ is $l_i$ restricted in $\overline{X}$.

Since $i(X)$ is a subspace of $\overline{X}$, $i:X \to \overline{X}$ is continuous. Thus, the same set function $\overline{i}: X \to \mathcal{D}\overline{X}$ is continuous. We have $em \circ \overline{i} = i$, where $em: (\mathcal{D}\overline{X},\{k_i\}_I) \to (Y,\{l_i\}_I)$ is the embedding map. Now, we need only to show that $em$ is a morphism in ${\bf DTop}(\Sigma,\mathcal{E})$.  The continuity is easily seen. For commutativity, we have

 $$\quad em(k_i(a_1,a_2,\dots,a_{n_i})) 
   = em(l_i(a_1,a_2,\dots,a_{n_i}))          
   = l_i(em(a_1),em(a_2),\dots,em(a_{n_i})).  $$
Thus, $U$ satisfies the solution set condition.
$\Box$
\vskip 3mm

By Theorem \ref{adjoint}, Lemma \ref{left1} and Lemma \ref{left2}, we get the main result of this section.
\begin{thm} \rm 
For every finitary signature $\Sigma$ and a set $\mathcal{E}$ of inequalities, the forgetful functor $U : {\bf DTop}(\Sigma,\mathcal{E}) \to {\bf DTop} $ has a left adjoint. Equivalently, For each directed space $X$, the free algebra over $X$ with respect to $\Sigma$ and $\mathcal{E}$ exists.
\end{thm}

\section{Continuous spaces as retracts of algebraic spaces}

In classical domain theory, it is shown that continuous domains are exactly the Scott retracts of algebraic domains \cite{CON03}. Besides, a dcpo $E$ is continuous iff the supremum map from $ID(E)$ to $E$ has a lower adjoint, which is the so-called way-below map \cite{CON03,KO95}.
We will construct a canonical algebraic space for each $T_0$ space and show that a $T_0$ space $X$ is continuous if and only if it is a continuous retract of an algebraic space if and only if the supremum map $\bigvee : I_T(X) \to X$ has a continuous lower adjoint $\DDa$.

\begin{defn} \rm Let $X$ be a $T_0$ space. A directed subset $D\subseteq X$ is called an ideal net if there exists  $x\in X$ such that:
\begin{enumerate}
\item[(1)] $D\rightarrow x$,
\item[(2)] $\forall d\in D$, $d\sqsubseteq x$.
\end{enumerate}
\end{defn}

Obviously, if $X$ is a poset endowed with the Scott topology, then every directed subset with a existing supremum is an ideal net.

\begin{lem} \rm If $D\subseteq X$ is an ideal net of a  $T_0$ space  $X$, then $\bigvee D$ exists and $D\rightarrow \bigvee D$.
\end{lem}

\noindent{\bf Proof.} By Lemma \ref{topo ideal}.

\begin{defn} \rm Let $X$ be a $T_0$ space. A subset $A$ of $X$ is called a topological ideal if there is an ideal net $D\subseteq X$ such that $A=\da D$. \end{defn}

So a topological ideal $A$ is a directed lower subset such that it has a supremum (denoted by $\bigvee A$) and converges to its supremum. Particularly, we set ${\rm I}_T(X)$ to be the set of all topological ideals of a $T_0$ space $X$.  Obviously, we have  $\{\da x: x\in X\}\subseteq {\rm I}_T(X)$. Set
$$\Omega({\rm I}_T(X))=\{{\mathcal U}\subseteq {\rm I}_T(X):  A\in {\mathcal U}\Leftrightarrow  \exists x\in X, \ \da x\in {\mathcal U} \ \&\ \da x\subseteq A\}.$$
Particularly, for any $x\in X$, we set ${\mathcal U}_x=\{A\in {\rm I}_T(X): x\in A\}$, then  ${\mathcal U}_x\in \Omega({\rm I}_T(X))$. We also denote $\rm I_T(X)$ to be $({\rm I}_T(X),\Omega({\rm I}_T(X)))$ if there is no confusion.

\begin{lem} \rm 
Let $X$ be a $T_0$ space. Then $({\rm I}_T(X),\Omega({\rm I}_T(X)))$ is an algebraic space such that:
\begin{enumerate}
\item[(1)] Its specialization order is equal to the set-theoretical inclusion;
\item[(2)] $K({\rm I}_T(X))=\{\da x: x\in X\}$.
\end{enumerate}
\end{lem}

\noindent{\bf Proof.}  Obviously, $\emptyset, {\rm I}_T(X) \in\Omega({\rm I}_T(X))$ and $\Omega({\rm I}_T(X))$ is closed under unions.  For ${\mathcal U},{\mathcal V}\in \Omega({\rm I}_T(X))$ and  $A\in {\mathcal U}\cap {\mathcal V}$,  there are $x,y\in X$ such that $\da x\in \mathcal U, \ \da y\in \mathcal V$ and $x,y\in A$. Since $A$ is directed, there exists $a\in A$ with $x,y\sqsubseteq a$. Thus $\da a\in \mathcal U\cap \mathcal V$ and $\da a\subseteq A$, i.e., $\mathcal U\cap \mathcal V\in \Omega({\rm I}_T(X))$. Hence, $\Omega({\rm I}_T(X))$ is a topology. If  $A\not\subseteq B$ for $A,B\in {\rm I}_T(X)$, then there is $a\in A$ together with $a\not\in B$; thus $A\in {\mathcal U}_x$ and $B\not\in{\mathcal U}_x$. If  $A \subseteq B$ for $A,B\in {\rm I}_T(X)$ and $A \in \mathcal U \in \Omega({\rm I}_T(X))$, then there exists some $\da x$ such that $\da x \subseteq A , \da x \in \mathcal U$; we have $\da x \subseteq B$ and thus $B \in \mathcal{U}$. Hence, $({\rm I}_T(X),\Omega({\rm I}_T(X)))$ is  $T_0$ and its specialization order is equal to the set-theoretical inclusion. 

Let $x \in X$, $\mathcal B$ be a directed subset of ${\rm I}_T(X)$ and $\mathcal{B} \to  \da x$. Since $\mathcal U_{x}$ is open, there must exists some $B \in \mathcal B $ such that $B \in \mathcal U_{x}$, i.e., $\da x \subseteq  B $. Hence, we have that $\da x$ is a compact element of $({\rm I}_T(X), \Omega({\rm I}_T(X))))$. For any $A\in {\rm I}_T(X)$, $\mathcal{A}=\{\da x:x\in A\}$ is directed. For $\mathcal{U}\in \Omega({\rm I}_T(X))$ with $A\in\mathcal{U}$, there exists $a\in A$ such that $\da a\in \mathcal{U}$. Thus,  $\mathcal{A}$ converges to $A$ in ${\rm I}_T(X)$. It follows that $({\rm I}_T(X),\Omega({\rm I}_T(X)))$ is an algebraic space and  $K({\rm I}_T(X))=\{\da x: x\in X\}$. $\Box$

\vskip 3mm

By definition, every topological ideal is directed and converges to its supremum. Hence, there exists a surjective map from ${\rm I}_T(X)$ onto $X$ as follows:
$$\bigvee: {\rm I}_T(X)\longrightarrow X, \ \ \forall A\in {\rm I}_T(X), \ A\mapsto \bigvee A.$$

\begin{prop} \rm
For a $T_0$ space $X$,  $\bigvee: {\rm I}_T(X)\longrightarrow X$ is a continuous surjective map.
\end{prop}

\noindent{\bf Proof.} Obviously, $\bigvee: {\rm I}_T(X)\longrightarrow X$ is surjective. For any $U\in \mathcal{O}(X)$ and for $A\in\bigvee^{-1}(U)$,
we have $ \bigvee A\in U$. Note that since $A$ is a topological ideal, we have that $A=\da A$ and $A\rightarrow \bigvee A$ in $X$. There exists $a\in A$ such that $a\in U$. Thus, $\da a\in \bigvee^{-1}(U)$ and $\da a\subseteq A$. Hence,  $\bigvee^{-1}(U)\in \Omega({\rm I}_T(X))) $, i.e., $\bigvee: {\rm I}_T(X)\longrightarrow X$ is continuous. $\Box$

\vskip 3mm

Hence, every $T_0$ space is a continuous image of some algebraic space. Particularly, when $X$ is a dcpo endowed with the Scott topology $\sigma(X)$,  ${\rm I}_T(X)$ is exactly the ideal completion ${\rm ID}(X)$ of  $X$ (see \cite{CON03}).

\begin{prop} \rm
Let $X$ be a dcpo endowed with the Scott topology $\sigma(X)$. Then ${\rm I}_T(X)= {\rm ID}(X)$ and $\Omega({\rm I}_T(X)) = \sigma({\rm I}_T(X))$, where ${\rm ID}(X)$ is the set of all ideal of $X$ endowed with the inclusion order.
\end{prop}
\noindent{\bf Proof.} Given any topological ideal $D$, it is an ideal by definition. Thus, $D \in {\rm ID}(X)$. Conversely, suppose that $A$ is a ideal of $X$. Since $X$ is a dcpo endowed with the Scott topology and $A$ is a directed subset, we know that $A \to \bigvee A$. Hence, $A$ is a topological ideal. Therefore, ${\rm I}_T(X)= {\rm ID}(X)$. 

Since ${\rm I}_T(X)$ is an algebraic space and $K({\rm I}_T(X)) = \{\da x : x \in X\}$,  $\{\big\uparrow_{{\rm I}_T(X)} (\da x)\}_{x \in X}$ form a base of the topology $\Omega({\rm I}_T(X))$ by Proposition \ref{basealg}. We know that $\{\big\uparrow_{{\rm ID}(X)} (\da x)\}_{x \in X}$ is also a base of $\sigma({\rm ID}(X))$ (see \cite{CON03}). Thus, $\Omega({\rm I}_T(X)) = \sigma({\rm I}_T(X))$.
$\Box$

\vskip 2mm
The following is the definition of adjunctions between posets.

\begin{defn} \rm
\cite{CON03}
Let $f:P\longrightarrow Q$ and  $g: Q\longrightarrow P$ be a pair of monotone maps between posets $P$ and $Q$. We say that $(f,g)$ is a pair of adjunctions provided the following holds:
$$\forall x\in P,y\in Q, \ f(x)\leq y\Leftrightarrow x\leq g(y).$$
In this case, $f$ is called  the lower adjoint of  $g$ and  $g$ is called the upper adjoint of $f$.
\end{defn}

\begin{thm} \rm \label{retract}
Let $X$ be a  $T_0$ space. Then following conditions are equivalent to each other:
\begin{enumerate}
\item[(1)] $X$ is a continuous space;
\item[(2)] $\bigvee: {\rm I}_T(X)\longrightarrow X$ has a continuous lower adjoint;
\item[(3)] $X$ is a continuous retract of some algebraic space.
\end{enumerate}
\end{thm}

\noindent{\bf Proof.} $(1)\Rightarrow (2)$. Let  $X$ is a continuous space. Then   $\dda x$ is directed and converges to $x$ for any $x\in X$ by Theorem. Note that since $\dda x=\da(\dda x)$, we have that  $\dda x$ is a topological ideal. Define a mapping  $\dda: X\longrightarrow {\rm I}_T(X)$ as follows:  $$\forall x\in X, \ \dda (x)=\dda x.$$
Then $\dda: X\longrightarrow {\rm I}_T(X)$ is well defined and monotone. For any $a\in X$ and for  $x\in \dda^{-1}({\mathcal U}_a)$, we have  $a\in \dda x$. Thus, there is some $y\in X$ such that  $a\ll y\ll x$. Set  $U=\UUa y$. Then  $U$ is an open neighborhood of  $x$ with  $U\subseteq \dda^{-1}({\mathcal U}_a)$. Hence, $\dda^{-1}({\mathcal U}_a)$ is an open subset of  $X$. Note that since $\{{\mathcal U}_a:a\in X\}$ is a base of  $\Omega({\rm I}_T(X))$, we have that  $\DDa: X\longrightarrow {\rm I}_T(X)$ is a continuous map. For any $x\in X$ and $A\in {\rm I}_T(X)$, $x\leq \bigvee A$ if and only if  $\DDa x\subseteq A$. Therefore,  $\DDa: X\longrightarrow {\rm I}_T(X)$ is the lower adjoint of  $\bigvee: {\rm I}_T(X)\longrightarrow X$ by definition.

$(2)\Rightarrow (3)$. Let $f:X\longrightarrow {\rm I}_T(X)$ be a continuous lower adjoint of $\bigvee: {\rm I}_T(X)\longrightarrow X$. Then  $\bigvee\circ f=id_X$ for that $\bigvee$ is surjective. Hence, $X$ is a continuous retract of the algebraic space ${\rm I}_T(X)$.

$(3)\Rightarrow (1)$. Let $Y$ be an algebraic space such that $X$ is a continuous retract of $Y$. Then there exist two continuous maps  $f:Y\longrightarrow X$ and  $g:X\longrightarrow Y$ with  $f\circ g=id_X$. For any $x\in X$ and $k\in K(Y)$, if  $k\leq g(x)$, then $f(k)\leq f(g(x))=x$. Let  $(x_j)_{j\in J}$ be a net of  $X$ with $(x_j)\rightarrow x$. Since  $g$ is continuous, we have  $(g(x_j))\rightarrow g(x)$. As  $k$ is compact and  $k\leq g(x)$, there exists  $j_0\in J$ such that  $k\leq g(x_{j_0})$. Hence, $f(k)\leq f(g(x_{j_0}))=x_{j_0}$, i.e., $f(k)\ll x$ by Definition. Note that since  $Y$ is an algebraic space, $\da g(x)\cap K(Y)$ is directed and converges to $g(x)$. Thus, $\{f(z):z\in \da g(x) \cap K(Y)\}\subseteq \DDa x$ and $(f(z))_{z \in \downarrow g(x) \cap K(Y)}\rightarrow f(g(x))=x$. Therefore, $X$ is a continuous space. $\Box$

\vskip 3mm

Since the lower adjoint of a map is unique, we have the following statement.
\begin{coro} \rm
A $T_0$ space $X$ is continuous iff the way-below map $\DDa : X \to I_T(X)$ is a well defined continuous map.
\end{coro}

\section{Algebraic spaces with way-below preserving maps}

In this section, we first give an equivalent order-theoretical description for algebraic spaces. Then we show that the category of all algebraic spaces with way-below preserving maps as morphisms is cartesian closed. Given any poset $A$, denote $ID(A)$ to be the set of all ideals of $A$ and denote $PI(A)$ to be the set of all principle ideals of $A$.

\begin{defn} \rm 
Let $X,Y$ be two algebraic spaces and $f:X \to Y$ be a continuous map. $f$ is said to be way-below preserving if for any $x \ll y$ in $X$, $f(x) \ll f(y)$ in $Y$.  
\end{defn}

\begin{lem} \rm 
Let $X,Y$ be two algebraic spaces. A continuous map $f:X \to Y$ is way-below preserving if and only if $f(K(X)) \subseteq K(Y)$.
\end{lem}
\noindent{\bf Proof.} 
Suppose that $f$ is way-below preserving. For any $x \in K(X)$, since $x \ll x$, $f(x) \ll f(x)$, i.e., $f(x) \in K(Y)$. 
For the converse, suppose that $f(K(X)) \subseteq K(Y)$. Since $X$ is algebraic, given any $x \in X$, $f(x) = \sup f(\da x \cap K(X))$. Given any $x \ll y$ in $X$, since $X$ is algebraic, there exists some $z \in K(X)$ such that $x \leq  z \ll z \leq y$. Then $f(x) \leq f(z) \ll f(z) \leq f(y)$, i.e., $f(x) \ll f(y)$.  $f$ is way-below preserving. $\Box$

\vskip 3mm

It is easy to check that all algebraic spaces with way-below preserving continuous maps as morphisms form a category, denoted by
$\bf ALG_\ll$. Now we define the notion of a b-poset.

\begin{defn} \rm 
We say that $(A,I(A))$ is a b-poset if $A$ is a poset and $PI(A) \subseteq I(A) \subseteq ID(A)$. Let $(B,I(B))$ be a b-poset. A monotone map $f:A \to B$  is called a b-map from $(A,I(A))$ to $(B,I(B))$ if $\da f(D) \in I(B)$ for any $D \in I(A)$.
\end{defn}

Given any b-maps $f:(A,I(A)) \to (B,I(B))$ and $g:(B,I(B)) \to (C,I(C))$, define $g\circ f :(A,I(A)) \to(C,I(C))$ with $g\circ f(a) = g(f(a))$ for any $a \in A$. Then $g\circ f$ is also a b-map. Thus, we can see that all b-posets together with  b-maps as morphisms form a category, denoted by $\bf b$-$\bf Poset$. 

Let $X$ be an algebraic space and $K(X)$ be the set of all compact elements of $X$. We denote $\mathcal{G}X$ to be a b-poset $(A,I(A))$ as follows: $A = K(X)$; $I(A) = \{\da\! x\cap K(X) : x \in X \}$. Let $Y$ be an algebraic space and $\mathcal{G}Y = (B,I(B))$.
For any way-below preserving map $f: X \to Y$, define $\mathcal{G}f : (A,I(A)) \to ((B,I(B)))$ to be a b-map as follows: $\mathcal{G}f(x) = f(x), \forall x \in A = K(X)$.

\begin{lem} \label{G functor}
$\mathcal{G}$ is a functor from $ \bf ALG_{\ll}$ to $\bf b$-$\bf Poset$.
\end{lem}
\noindent{\bf Proof.} Let $X$ be an algebraic space.
Then, for any $x \in X$, $\da\! x \cap K(X)$ is an ideal of $K(X)$. For any $x \in K(X)$, $\da\! x \cap K(X) = \{y \in K(X): y \leq x \}$ is a principle ideal of $K(X)$.
Thus, $\mathcal{G}X$ is indeed a b-poset. Suppose that $X,Y$ are two algebraic spaces and $\mathcal{G}X = (A,I(A)), \mathcal{G}Y = (B,I(B))$. For any way-below preserving continuous map $f:X \to Y$ and $x \in K(X)$, $\mathcal{G}f(x) = f(x) \in K(Y) = B$. Obviously, $\mathcal{G}f$ is monotone.  For any $D =\ \da\! x \cap K(X) \in I(A)$, we have $\da \mathcal{G}f(D) =\ \da\! f(D) = f(x) \cap K(Y) \in I(B)$. Thus, $\mathcal{G}f$ is a b-map between $\mathcal{G}X$ and $\mathcal{G}Y$. It is easy to see that any identity morphism is preserved by $\mathcal{G}$ and $\mathcal{G}(g \circ f) = \mathcal{G}(g) \circ \mathcal{G}(f)$. Therefore, $\mathcal{G}$ is a functor.
$\Box$

\vskip 3mm

Conversely, given any b-poset $(A,I(A))$, we define $\mathcal{H}(A,I(A))$ to be an algebraic space $X$ as follows. The carrier set of $X$ is $I(A)$. For any $a \in A$, denote $\mathcal{B}(a) = \{D \in I(A): a \in D\}$. The topology of $X$ is generated by taking $\mathcal{B}_A = \{\mathcal{B}(a): a \in A\}$ as a base. Given any b-map $f:(A,I(A)) \to (B,I(B))$, $\mathcal{H}f$ is defined by $\mathcal{H}f(D) = \da f(D), \forall D \in I(A) =X$. 
\begin{lem} \label{H functor}
$\mathcal{H}$ is a functor from $\bf b$-$\bf Poset$ to $\bf ALG_{\ll}$
\end{lem}
\noindent{\bf Proof.} Let $(A,I(A))$ be any b-poset and $X = \mathcal{H}(A,I(A))$. First, we check that $\mathcal{B}_A$ form a base of a topology on $I(A)$. Given any $ \mathcal{B}(a)$ and $\mathcal{B}(b)$ in $\mathcal{B}_A$, assume that $D \in \mathcal{B}(a) \cap \mathcal{B}(b)$, then $a,b \in D$. Since $D$ is an ideal of $A$, there exists a $c \in D$ such that $a,b \leq c$ in $A$. Thus, $D \in \mathcal{B}(c) \subseteq \mathcal{B}(a) \cap \mathcal{B}(b) $, i.e., $\mathcal{B}_A$ is a base. It is also easy to check that the specialization order of the topology generated by $\mathcal{B}_A$ is equal to the set inclusion order on $I(A)$. Then for each $D \in I(A)$, $\{\da a : a \in D\}$ is a directed subset of $I(A)$ and converges to $D$ relative to $\mathcal{B}_A$. Moreover, for any $a \in A$, since $\mathcal{B}(a) = \{D \in I(A): a \in D\} = \big\ua_{I(A)}( \da a)$ is an open subset, $\da a$ is a $n$-compact element of $X$. Thus, $X$ is an algebraic space with $K(X) = PI(A)$.

Let $Y = \mathcal{H}(B,I(B))$ and $f:(A,I(A)) \to (B,I(B))$ be a b-map. Then $\mathcal{H}f$ is a well-defined map from $X $ to $Y$. Now, we show that it is continuous and way-below preserving. Given any $\mathcal{B}(b) \in \mathcal{B}_B$, $\mathcal{H}f^{-1}(\mathcal{B}(b)) = \{D \in I(A): b \in\ \da f(D)\} = \{D \in I(A): \exists d \in D. b \leq f(d)\} = \bigcup_{ b \leq f(d)} \mathcal{B}(d)$ is open in $X$. Therefore, $\mathcal{H}f$ is continuous. Since the set of compact elements of $X$ and $Y$ are respectively the principle ideals of $A$ and $B$, $\mathcal{H}f(K(X)) \subseteq K(Y)$. Thus, $\mathcal{H}f$ is way-below preserving. It is easy to check that $\mathcal{H}$ preserves every identity morphism and the composition of morphisms.
$\Box$

\vskip 3mm

By definition, we have that $$\mathcal{H}\mathcal{G}X = \mathcal{H}(K(X),\{\da x \cap K(X)\}_{x \in X}) = (\{\da x \cap K(X)\}_{x \in X}, \mathcal{B}_{K(X)})$$ is isomorphic to  $X$, and $$\mathcal{G}\mathcal{H}(A,I(A)) = \mathcal{G}(I(A),\mathcal{B}_{A}) = (\{\da a : a\in A\},\{\da x \cap \{\da a : a\in A\} \}_{ x \in I(A)}) $$ is isomorphic to  $ (A,I(A))$. Thus, there is a one-to-one correspondence between algebraic spaces (i.e., b-spaces) and b-posets up to isomorphism. Also, a map is a b-map between b-posets iff it is a way-below preserving map for the corresponding b-spaces. Thus, we have the following statement.
\begin{prop}
 $\bf b$-$\bf Poset$ is equivalent to $\bf ALG_{\ll}$
\end{prop}
\noindent{\bf Proof.} Easily seen from Lemma \ref{G functor} and Lemma \ref{H functor} and the above discussion. $\Box$

\vskip 3mm


The following example shows that there exist two algebraic spaces which are not homeomorphic, while the posets of their compact elements are order isomorphic. 

\begin{exmp} \rm 
Let $X = \mathbb{N}^{\top}$ be the poset of natural numbers adding a top element $\top$. Let $Y = \mathbb{N} \cup \{\omega_1,\omega_2\}$, where $\omega_1$ is bigger than all natural numbers and $\omega_2$ is bigger than $\omega_1$. Considering $(X,A(X))$ and $(Y,\sigma(Y))$, they are both algebraic spaces. The sets of compact elements of them are $\mathbb{N}^\top$ and $\mathbb{N} \cup \{\omega_2\}$ respectively, which are order isomorphic.
\end{exmp}

\vskip 3mm

We say that two algebraic spaces $X,Y$ are basis isomorphic with basis $A$ if $K(X)$ and $K(Y)$ are both order isomorphic to a poset $A$. By viewing algebraic spaces as b-posets, all algebraic spaces with a basis isomorphic to $A$ are just all the b-posets $(A,I(A))$ with $PI(A) \subseteq I(A) \subseteq ID(A)$ up to isomorphism.

The following statement is the main result of this section.

\begin{thm}
$\bf b$-$\bf Poset$ is cartesian closed.
\end{thm}

\noindent{\bf Proof.} Obviously, $\{*,\{*\}\}$ is the terminal object. Given any two objects $(A,I(A))$ and $(B,I(B))$, denote $(A,I(A)) \times (B,I(B))$ to be $(A \times B, I(A \times B))$, where $D \in I(A \times B)$ iff $\pi_A (D) \in I(A),\pi_B (D) \in I(B)$ and $D = \pi_A (D) \times \pi_B (D)$, where $\pi_A,\pi_B$ are the projections from $A \times B$ to $A$ and $B$ respectively. Then $I(A \times B)$ is a set of ideals of $A \times B$ and contains all principle ideals of $A \times B$. Thus, $(A,I(A)) \times (B,I(B))$ is a well-defined b-poset. We show that it is the categorical product.

Given any b-poset $(C,I(C))$ and b-maps $f_1,f_2$ from $(C,I(C))$ to $(A,I(A))$ and $(B,I(B))$ respectively. It is easy to check that the projections $$\pi_A:(A \times B, I(A\times B)) \to (A,I(A)) , \pi_A((a,b)) = a, \forall(a,b) \in A \times B$$ and
$$\pi_B:(A \times B, I(A \times B)) \to (B,I(B)) , \pi_B((a,b)) = b, \forall(a,b) \in A \times B$$ 
are both b-maps, and $$\langle f_1, f_2 \rangle:(C,I(C)) \to (A \times B, I(A \times B)),$$
where $$\langle f_1, f_2 \rangle(c) = (f_1(c),f_2(c)),\forall c\in C,$$ 
is a b-map such that $\pi_A \circ \langle f_1, f_2 \rangle = f_1$ and  $\pi_B \circ \langle f_1, f_2 \rangle = f_2$. 
\vskip 2mm
Let $B^A$ be the poset (the pointwise order) of all b-maps from $(A,I(A))$ to $(B,I(B))$. 
Define the evaluation map $ev : B^A \times A \to B$ with $ev(f,a) = f(a)$.
Denote $I(B^A)$ to be the largest set of ideals of $B^A$ such that the evaluation map $ev$ is a well-defined b-map from $(B^A,I(B^A)) \times (A,I(A))$ to $(B,I(B))$.
Since for any principle ideals $\da f$ of $B^A$ and for any $E \in I(A)$, $\da ev(\da f \times E) =\ \da \{f(e): e \in E \} \in I(B)$, then $I(B^A)$ contains all principle ideals of $B^A$. $(B^A,I(B^A))$ is a well-defined b-poset.

We claim that $(B^A,I(B^A))$ is the exponential object for $(A,I(A))$ and $(B,I(B))$. 
Given any b-map $f: (C,I(C)) \times (A,I(A)) \to (B,I(B))$, let $\overline{f}: C \to B^A$ be a map as follows: for any $c \in C$, $\overline{f}(c) = f_{c}$, where $\forall a \in A, f_c(a) = f(c,a)$. Since for any $D \in I(A)$, $ \da f_c(D) =\ \da \{f(c,d) : d\in D \} =\ \da \{f(a,d): a \in \da c, d \in D\} =\ \da f(\da c \times D) \in I(B)$ , then $\overline{f}$ is a well defined map from $C$ to $B^A$ and it is monotone. For any $D \in I(C)$, $\overline{f}(D) = \da_{B^A} \{f_d : d\in D\}$ is an ideal of $B^A$. Given any $F\in I(A)$, $ ev(\overline{f}(D),F) =\ \da\! \{f(d,a) : (d,a) \in D\times F \} =\ \da f(D \times F) \in I(B)$, i.e., $\overline{f}(D) \in I(B^A) $. Thus, $\overline{f}$ is a b-map from $(C,I(C))$ to $(B^A,I(B^A))$.

It is easy to check that for any b-map $f:C \times A \to B$, $f = ev \circ (\overline{f} \times  id_A) $, where $ \overline{f} \times  id_B (c,a)= (\overline{f}(c),a), \forall (c,a) \in C\times A$. 
Suppose that $g:C \to B^A$ satisfies $ f = ev \circ (g \times  id_A)$. For any $(c,a) \in C \times A$, $ev\circ(g \times id_A)(c,a) = ev(g(c),a) = g(c)(a) = f(c,a)$. Therefore, $g(c) = f_c = \overline{f}(c)$ for any $c \in C$, $g = \overline{f}$. It follows that $(B^A,I(B^A))$ is the exponential object.
 $\Box$
\vskip 3mm

\begin{cor}
$\bf ALG_\ll$ is cartesian closed.
\end{cor}

Denote $\bf ALGD_\ll$ the category of all algebraic domains with way-below preserving maps as morphisms and $\bf Poset$ the category of all posets with monotone maps as morphisms.
It is well known that $\bf ALGD_\ll$ is equivalent to $\bf Poset$ \cite{YAN2004} and then is cartesian closed. In fact, they can both be viewed as subcategories of $\bf ALG_\ll$. Given any algebraic domains $P$, it can be viewed as an algebraic space $(P,\sigma(P))$ and it can also be viewed as a b-poset $(K(P),ID(K(P)))$. A way-below preserving map between two algebraic domains also corresponds to a b-map between the corresponding b-posets.

Similarly, a poset $P$ can be viewed as a b-poset $(P,PI(P))$. A monotone map $f: P \to Q$ corresponds to the b-map $f:(P,PI(P)) \to (Q,PI(Q))$.  For convenience, we denote $\bf ALGD^*_\ll$ and $\bf Poset^*$ to be the  subcategories of $\bf ALG_\ll$ corresponding to $\bf ALGD_\ll$ and $\bf Poset$ respectively, i.e., 
$\bf ALGD^*_\ll$ is the full subcategory of $\bf ALG_\ll$, for which objects are of the form $(A,ID(A))$ and $\bf Poset^*$ is the full subcategory of $\bf ALG_\ll$ for which objects are of the form $(A,PI(A)))$. Then, we have the following statement.

\begin{thm} \rm 
$\bf ALGD^*_\ll$ and $\bf Poset^*$ are both  cartesian closed subcategories of $\bf ALG_\ll$. Moreover, the  exponential objects in $\bf ALGD^*_\ll$ agree with those in $\bf ALG_\ll$.
\end{thm}
\noindent{\bf Proof.}
By the above discussing we know that $\bf ALGD^*_\ll$ and $\bf Poset^*$ are cartesian closed full subcategories of $\bf ALG_\ll$. Given any two b-posets $(A,ID(A))$ and $(B,ID(B))$,  their exponential object in $\bf ALG_\ll$ is $(B^A,I(B^A))$, where $I(B^A)$ is the largest set of ideals of $B^A$ that makes $ev:(B^A,I(B^A)) \times (A,ID(A)) \to (B,ID(B))$ to be a well-defined b-map. Since given any ideal $D$ of $B^A$ and ideal $E$ of $B$, $ev(D\times E) =\ \da \{f(e) : (f,e) \in D \times E\}$ is a ideal of $B$, we have $I(B^A) = ID(B^A)$. Thus, $(B^A,I(B^A))$ is also the exponential object of $(A,ID(A))$ and $(B,ID(B))$ in $\bf ALGD^*_\ll$. In other words, the exponential object in $\bf ALGD^*_\ll$ coincide with those in $\bf ALG_\ll$. $\Box$

\vskip 3mm

However, the exponential objects in $\bf Poset^*$ may not agree with those in $\bf ALG_\ll$. We give an example as follows. 

\begin{exmp} \rm
Consider the poset of natural numbers $\mathbb{N}$ and the  poset $B = \{\bot,\top\}$ with $\bot \leq \top$. In $\bf Poset^*$, the exponential object of $(\mathbb{N},PI(\mathbb{N}))$ and $(B,PI(B))$ is $(\mathbb{N}^B,PI(\mathbb{N}^B))$, where $\mathbb{N}^B$ is the set of all monotone maps from $\mathbb{N}$ to $B$ with the pointwise order. However, given any $D \in ID(\mathbb{N}^B)$ and $E \in ID(\mathbb{N})$, $ev(D,E) =\ \da \{f(n): (f,n) \in D \times E \}$ is always a principle ideal of $B$. Thus, the exponential object of $\mathbb{N}$ and $B$ in $\bf ALG_\ll$ is $(\mathbb{N}^B,ID(\mathbb{N}^B))$, which is different from $(\mathbb{N}^B,PI(\mathbb{N}^B))$.
\end{exmp}



Finally, we show that $\bf ALGD^*_\ll$ is a reflective subcategory of $\bf ALG_\ll$. And for any algebraic space $X$, the sobrification of $X$ is just the reflection of $X$ over $\bf ALGD_\ll$.

\begin{thm}
$\bf ALGD^*_\ll$ is a reflective subcategory of $\bf ALG_\ll$.
\end{thm}
\noindent{\bf Proof.}
Given any object $(A,I(A))$ of $\bf ALG_\ll$, let $i:A \to A$ be the identity map. Then $i$ is also a b-map from  $(A,I(A))$ to $(A,ID(A))$. For any object $(B,ID(B))$ of $\bf ALGD^*_\ll$ and any b-map $f:(A,I(A)) \to  (B,ID(B))$, let $\overline{f}:(A,ID(A)) \to (B,ID(B))$ with $\overline{f}(x) = f(x),\forall x\in A$. Then $\overline{f}$ is a b-map. It is easily seen that $f = \overline{f} \circ i$ and such $\overline{f}$ is unique. Thus, $\bf ALGD^*_\ll$ is a reflective subcategory of $\bf ALG_\ll$.

\begin{thm}\rm
For any algebraic space $X = \mathcal{H}(A,I(A))$, the sobrification of $X$ is equal to $\mathcal{H}(A,ID(A))$ up to homeomorphism.
\end{thm}
\noindent{\bf Proof.}
Since $\mathcal{H}(A,ID(A))$ is an algebraic domain endowed with the Scott topology,  it is a sober space. Since the basis of the topology of $(A,ID(A))$ is $\mathcal{B}(A)$ and is equal to that of $(A,I(A))$, their topology are order isomorphic. Thus, $\mathcal{H}(A,ID(A)))$ is equal to the sobrification of $X$ up to homeomorphism.
\vskip 3mm

The above theorem also tell us that for all basis isomorphic algebraic spaces, up to homeomorphism, their sobrification are the same and are equal to the algebraic domain which is basis isomorphic with them.

\vskip 3mm

\section{Category of continuous spaces}

In this section, we show that $\bf CON_\ll$, the category of continuous spaces with way-below preserving continuous maps as morphisms, form a cartesian closed category.

\begin{defn}\rm \cite{AJ94}
An abstract basis is given by a set $A$ together with a transitive relation $\prec$ on $A$ such that $$M \prec z \Longrightarrow \exists y \in A.\ M \prec y \prec z$$
holds for all elements $z$ and finite subsets $M$ of $A$.
\end{defn}

\begin{defn} \rm 
We call $(A,\leq,\prec)$ a normal abstract basis if $(A,\leq)$ is a poset, $(A,\prec)$ is an abstract basis and the following two conditions are satisfied: (1) $a \prec b$ implies $a \leq b$.
(2) if $a \leq b \prec c \leq d$, then $a \prec d$; (3) if $a \not\leq b$, then there exists a $c$ such that $c \prec a$ and $c \not \prec b$.
\end{defn}

\begin{lem}\rm \label{con abs}
Let $X$ be a continuous space. Then $(X,\sqsubseteq,\ll)$ is a normal abstract basis.
\end{lem}
\begin{proof}
Obviously.
\end{proof}

\begin{prop}\rm \label{abs topo}
Let $(A,\leq,\prec)$  be a normal abstract basis. We denote $\upp a$ to be $\{b: a \prec b\}$.
Then $\{\upp a: a \in A\}$ form a base of a topology, denoted by $\tau$, on $Y^X$. Moreover, $(A,\tau)$ is a continuous space. $a \leq b$ iff $a \sqsubseteq b$ relative to $(A,\tau)$.  $a \ll b $ in $(A,\tau)$ iff $a \prec b$.
\end{prop}
\noindent{\bf Proof.}  Given any $a,b \in A$ and $c \in\ \upp a\ \cap \upp b$, by Lemma \ref{prec int}, there exists a $d$ such that $a,b \prec d \prec c$. Thus, $d \in\ \upp a\ \cap \upp b$ and $c \in\ \upp d$. $\{\upp a: a \in A\}$ is  base of a topology. By definition, every open subset of $(A,\tau)$ is an upper set relative to $\leq$. Thus,  $a \leq b$ implies $a \sqsubseteq b$ relative to $(A,\tau)$. Assuming that $a \not \leq b$, then there exists a $c \in A$ such that $a \in\ \upp c$ and $b  \not \in\  \upp c$, then, $a \not\sqsubseteq b$. Therefore, $\leq$ is equal to $\sqsubseteq$.


Now, we show that $(A,\tau)$ is a continuous space. Given any $b \in\ \upp a$, since there exists a $c$ such that $c \in\ \upp a$ and $b \in\ \upp c$, then $b \in (\ua c)^\circ \subseteq\ \ua c \subseteq\  \upp a $. Therefore, $(A,\tau)$ is a continuous space.

Suppose $a \prec b$, then $b \in\ \upp a \subseteq (\ua a)^\circ$, i.e., $a \ll b$. Conversely, suppose $b \in (\ua a)^\circ$. Since $\{\upp a: a \in A\}$ is a base of the topology $\tau$, $\{c : c \prec b \}$ form a directed subset that converges to $b$ relative to $\tau$. Then there exists some $c \prec b$ such that $ c \in  (\ua a)^\circ$, then $a \leq c \prec b$, then $a \prec b$.
 $\Box$

\begin{cor}\rm 
There is a one-to-one correspondence between a normal abstract basis and a continuous space.
\end{cor}
\noindent{\bf Proof.} By Lemma \ref{con abs} and Proposition \ref{abs topo}
$\Box$

\begin{cor}\rm \label{con YX}
For any two continuous spaces $X$ and $Y$, the topological space $(Y^X,\tau)$ generated by $ \prec_0$ is a continuous space. The specialization order is the set inclusion order. The $\ll$ relation is equal to $\prec_0$. 
\end{cor}\rm
\noindent{\bf Proof.} Obviously, if $A \leq B \prec_0 C \leq D$ in $Y^X$, then $A \prec_0 D$. By Proposition \ref{abs topo}, we need only to show that if $A \not \leq B$, then there exists a $C$ such that $C \prec_0 A$ and $C \not\prec_0 B$. Assume $A =\ \da \{f_z:z \ll z_0\}$ for some continuous space $Z$, continuous map $f:Z \otimes X \to Y$ and $z_0 \in Z$. Since $A \not\leq B$, there exists some $z_1 \ll z_2 \ll z_0$ such that $f_{z_1} \in A$ and $f_{z_1} \not \in B$. Then $\da \{f_z:z\ll z_2\} \prec_0 A$ and $\da \{f_z:z\ll z_2\} \not \prec_0 B$. $\Box$

\begin{defn}\rm 
Let $(A,\leq_A,\prec_A)$ and $(B,\leq_B,\prec_B)$ be two normal abstract basis. A map $f:A \to B$ is said to be a normal map from $(A,\leq_A,\prec_A)$ to $(B,\leq_B,\prec_B)$ if it preserves $\leq$ and $\prec$ relation and for any $y \prec f(x)$, there is a $z \in X$ such that $z \prec x$ and $y \prec f(z)$.
\end{defn}

\begin{prop}
The normal maps and way-below preserving continuous maps are one-to-one corresponding.
\end{prop}
\begin{proof}
Let $(A,\leq_A,\prec_A)$ and $(B,\leq_B,\prec_B)$ be two normal abstract basis. And $(A,\tau_A),(B,\tau_B)$ be the corresponding continuous spaces. Assume that $f:(A,\leq_A,\prec_A) \to (B,\leq_B,\prec_B)$ is a normal map. Then it is a way-below preserving map from $(A,\tau_A)$ to $(B,\tau_B)$. For the continuity, consider $U = f^{-1}(\UUa y)$ for any $y \in B$. Supposing that $x \in U$, then $y \prec_B f(x)$ and there exists some $z \in X$ such that $z \prec x$ and $y \prec f(z)$, i.e., $x \in \UUa z$ and $z \in U$. Thus, $U$ is open and then $f$ is a way-below preserving continuous map from $(A,\tau_A)$ to  $(B,\tau_B)$.

Conversely, assume that $f$ is a way-below preserving continuous map from $(A,\tau_A)$ to  $(B,\tau_B)$. Obviously, $f$ preserves $\leq$ and $\prec$ relations. Supposing $y \prec_B f(x)$, then $y \ll_B f(x)$. Since $(A,\tau_A)$ is a continuous space and $f$ is continuous, then $\{f(z): z \ll_A x \}$ is a directed subset that converges to $f(x)$ relative to $\tau_B$. Since $(B,\tau_B)$ is a continuous space, there exists some $y_1 \in B$ such that $y \ll_B y_1 \ll_B f(x)$. It follows that there exists some $z_0 \ll_A x$ such that $y \ll_B y_1 \sqsubseteq f(z_0) $ relative to $\tau_B$, i.e., $y \prec_B f(z_0)$. Therefore, $f$ is a normal map from  $(A,\leq_A,\prec_A)$ to $(B,\leq_B,\prec_B)$.
\end{proof}

Given any two normal maps $f: (A,\leq_A,\prec_A) \to (B,\leq_B,\prec_B)$ and $g: (B,\leq_B,\prec_B) \to (C,\leq_C,\prec_C)$. Then $g \circ f:A \to C$ is a normal map from $(A,\leq_A,\prec_A)  \to (C,\leq_C,\prec_C)$.
Denote $\bf NAB$ the category of normal abstract basis with normal maps as morphisms. 

\begin{cor}
$\bf CON_\ll$ is equivalent to $\bf NAB$.
\end{cor}

In the following, we construct the exponential object in $\bf CON_\ll$ and show that $\bf CON_\ll$ is cartesian closed.

Let $X,Y,Z$ be continuous spaces and $f: Z \otimes X \to Y$ be a way-below preserving continuous map. We define a map $\overline{f}: Z \to ID([X \to Y])$ with $\overline{f}(z_0) =\ \da \{f_{z}: z\ll z_0\}$, where $ID([X \to Y])$ is the set of ideals of $[X \to Y]$, $[X \to Y]$ is the poset of all continuous maps from $X$ to $Y$ equipped with the pointwise order, and $f_z: X \to Y$ with $f_z(x) = f(z,x), \forall x \in X$. Denote $Y^X$ a subset of $ID([X \to Y])$ that $A \in Y^X$ iff there are a continuous space $Z^\prime$, a way-below preserving continuous map $g:Z^\prime \otimes X \to Y$ and a $z_1 \in Z^\prime$ such that $A = \overline{g}(z_1)$. Note that for continuous spaces $X,Y$, $X \otimes Y = X \times Y$ by Proposition \ref{core times}

\begin{defn}\rm {}
Given any two continuous spaces $X,Y$, we define a relation $\prec_0$ on $Y^X$ as follows: $A \prec_0 B$ iff there exists $h \in B$ such that $A \subseteq \da h \subseteq B$, where $\da h$ means $\{f \in [X \to Y]: f \leq h\}$.
\end{defn}

\begin{lem}\rm  \label{prec int}
For any two continuous spaces $X$ and $Y$, $(Y^X,\prec_0)$ is an abstract basis. 
\end{lem}
\noindent{\bf Proof.} Obviously, $\prec_0$ is transitive. Assume that $A_1,A_2 \prec_0 B$ in $Y^X$ and $B =\ \da \{f_{z_0}:z \ll z_0\}$ for some continuous space $Z$, way-below preserving continuous map $f:Z \otimes X \to Y$ and $z_0 \in Z$. By definition, there exists some $z_1,z_2 \ll z_0$ such that $A_1 \subseteq 
\da h_1 \subseteq \da f_{z_1} \subseteq B $ and $A_2 \subseteq 
\da h_2 \subseteq \da f_{z_2} \subseteq B $. There exists some $z_3,z_4$ such that $z_1,z_2 \ll z_3 \ll z_4 \ll z_0$ and it follows that $A_1,A_2 \subseteq \da f_{z_3}  \subseteq \da \{f_{z}:z \ll z_4\} \subseteq \da f_{z_4} \subseteq B$, i.e., $A_1,A_2 \prec_0\ \da \{f_{z}:z \ll z_4\} \prec_0 B$. $(Y^X,\prec_0)$ is an abstract basis. $\Box$

\begin{prop}\rm \label{over f}
Let $X,Y,Z$ be continuous spaces and $f: Z \otimes X \to Y$ be a way-below preserving continuous map. Then
$\overline{f}: Z \to Y^X$ is a way-below preserving continuous map, where $Y^X$ is endowed with topology $\tau$ generated by $\prec_0$.
\end{prop}
\noindent{\bf Proof.} Recall that $\overline{f}(z_0) =\ \da\! \{f_z:z\ll z_0\}$. If $z_1 \ll z_2$ in $Z$, there exists some $z_3$ such that $z_1 \ll z_3 \ll z_2$. It follows that  $\overline{f}(z_1) =\ \da \{f_z:z\ll z_1\} \subseteq\ \da f_{z_3} \subseteq\ \da \{f_z:z\ll z_2\} =\ \overline{f}(z_2)$, i.e., $\overline{f}(z_1) \ll \overline{f}(z_2)$ in $(Y^X,\tau)$. Therefore, $\overline{f}$ is way-below preserving.

For the continuity of $\overline{f}$, we need only to show that $\overline{f}^{-1}(\UUa A)$ is open for any $A \in Y^X$, where $\UUa A = \{B \in Y^X: A \ll B \text{ relative to } \tau \}$. Obviously, $\overline{f}$ is monotone. Suppose that $z_0 \in \overline{f}^{-1}(\UUa A)$, i.e., $A \ll\ \da \{f_z:z \ll z_0\}$. Then, there exists some $z_1 \ll z_0$ such that $A \subseteq\ \da f_{z_1}$. Since there exists some $z_2$ such that $z_1 \ll z_2  \ll z_0$, then $A \subseteq\ \da f_{z_1} \subseteq \overline{f}(z_2) \subseteq \overline{f}(z_0)$.  Thus, $z_2 \in \overline{f}^{-1}(\UUa A)$ and $z_0 \in \UUa z_2$. $\overline{f}^{-1}(\UUa A)$ is open. $\Box$

\vskip 3mm

We define the evaluation map $ev: Y^X \otimes X \to Y$ as $ev(A,x) = \sup \{g(x):g \in A\}$.

\begin{prop}\rm \label{def ev} The map $ev$ is a  well-defined way-below preserving continuous map.
\end{prop}
\noindent{\bf Proof.} For any $(A,x) \in Y^X \otimes X$, assume that $A =\ \da \{f_z:z \ll z_0\}$ for some $f,Z$ and $z_0$. Since $f:Z \otimes X \to Y$ is continuous, then $ev(A,x) = \sup \{f_z(x):z \ll z_0\} = \sup \{f(z,x): z \ll z_0\} = f(z_0,x)$ is well-defined and $\{f_z(x): z \ll z_0\} \to  \sup \{f_z(x): z \ll z_0\}$.

Suppose $(A,x) \ll (B,y)$ in $Y^X \otimes X$ and $B =\ \da \{f_z:z \ll z_0\}$. Then $A \ll B, x \ll y$ and there exists some $h$ such that  $A \subseteq \da h \subseteq B$. $ev(A,x) = \sup \{g(x):g \in A\} \leq h(x) \leq f_{z_1}(x) = f(z_1,x) \ll f(z^\prime,y) \leq \sup\{f_z(y):z \ll z_0\} = ev(B,y)$ for some $z_1 \ll z^\prime \ll z_0$. Therefore, $ev$ is way-below preserving.

For continuity of $ev$, we need only to show that $ev$ is separately continuous by Proposition \ref{separate}

(1) Consider $ev(A,\_ ): X \to Y$ for any $A \in Y^X$. Assume that $A =\ \da\{f_z:z \ll z_0\}$ for some continuous space $Z$, way-below preserving map $f:Z \otimes X \to Y$ and $z_0 \in Z$. 
Then $ev(A,x) = f(z_0,x)$. Thus, $ev(A,\_ )$ is continuous.


(2) Consider $ev(\_,x) : Y^X \to X$ for any $x \in X$. For any $y \in X$, let $U = ev(\_,x)^{-1}(\UUa y)$ and $A \in U$. Suppose $A = \overline{f}(z_0) =  \da \{f_z:z \ll z_0\}$ for some $Z,f,z_0$. Then there exists some $y^\prime$ such that $y \ll y^\prime \ll \sup\{f_z(x):z \ll z_0\}$. It follows that there exists some $z^\prime \ll z_0$ such that $y \ll f_{z^\prime}(x) = \sup\{f_{z}(x):z \ll z^\prime\}$. Then $\overline{f}(z^\prime) \in U$ and $\overline{f}(z^\prime) \ll \overline{f}(z_0) = A$. Therefore, $U$ is open and $ev(\_,x)$ is continuous. $\Box$

\begin{prop}\rm \label{ev}
$ev(\overline{f} \otimes id) = f$. For any $g$ such that $ev(g \otimes id) = f$, we have $g = \overline{f}$.
\end{prop}
\noindent{\bf Proof.} We have
\begin{align*}
  & ev(\overline{f} \otimes id)(z_0,x_0) \\
=\ & ev(\overline{f}(z_0),x_0) \\
=\ & ev(\da \{f_z:z \ll z_0\},x_0) \\
=\ & \sup\{f(z,x_0):z \ll z_0 \} \\
=\ & f(z_0,x_0).
\end{align*}
Suppose that $g:Z \to Y^X$ be a way-below preserving continuous map such that $ev(g \otimes id) = f$. We show $g = \overline{f}$.

Given any $z_0 \in Z$, $g(z_0) = \sup \{g(z):z \ll z_0\}$. For any $z_1 \ll z_0$, we have $\overline{f}(z_1) \ll \overline{f}(z_0)$, supposing that $\overline{f}(z_1) \subseteq \da h \subseteq \overline{f}(z_0)$. Then for any $x \in X$,
\begin{align*}
  & ev(g \otimes id)(z_1,x) \\
 =\ & ev(g(z_1),x) \\
 =\ & \sup\{k(x): k \in g(z_1) \}\\
 =\ & ev(\overline{f} \otimes id)((z_1,x) \\
=\ & \sup\{f_z(x):z \ll z_1 \} \\
\leq\ & h(x). 
\end{align*}
Thus, $g(z_1) \subseteq \da h \subseteq \overline{f}(z_0)$. It follows that $g(z_0) = \sup\{g(z):z \ll z_0\} \leq \overline{f}(z_0)$. Then $g \leq \overline{f}$. The same, we can show $\overline{f} \leq g$. Thus, $f = g$.


\begin{thm}
$\bf CON_\ll$ is a cartesian closed category.
\end{thm}
\noindent{\bf Proof.} By Corollary \ref{con YX}, Proposition \ref{over f}, Proposition \ref{def ev} and Proposition \ref{ev}.

\section{Free algebras over continuous spaces}

In this section, we use a categorical method, which was originated from \cite{KO95}, to prove that the carrier space of a free algebra over a continuous space with respect to ${\bf DTop}(\Sigma,\mathcal{E})$ is a continuous space. The result in the previous section that continuous spaces are retracts of algebraic spaces will play an important role in the proof.

\begin{lem} \rm 
Define an functor $T$ as follows: 

$$\forall X \in {\rm Ob}({\bf DTop}),\ T(X) = {\rm I}_{T}(X);$$
$$\forall f \in {\rm Mor}({\bf DTop}),\ \forall D \in T(X),\ T(f)(D) = \da f(D).$$
Then, $T$ is a endofunctor on ${\bf DTop}$.
\end{lem}

\noindent{\bf Proof.} First, we show that $T(f)$ is well defined and is a continuous map from ${\rm I}_{T}(X)$ to ${\rm I}_{T}(Y)$. Given any $D \in {\rm I}_{T}(X)$, we have $T(f)(D) = \da f(D) \to f(\sup D)$, since $D \to \sup D$. Thus, $f(\sup D)$ is the supremum of $\da f(D)$ by Lemma \ref{topo ideal}. $T(f)(D) \in {\rm I}_T(Y)$. 

Since ${\rm I}_{T}(Y)$ is a algebraic space with compact elements $\{\da y: y \in Y\}$, we need only to show that $T(f)^{-1} (\ua_{{\rm I}_{T}(Y)} (\da y))$ is open for any $y \in Y$. Indeed, we have that
\begin{align*}
& T(f)^{-1} (\ua_{I_{T}(Y)} (\da y)) \\
=\ & \{D \in I_{T}(X) : y \in \da f(D)\} \\
=\ & \{D \in I_{T}(X) : \exists d \in D.\ y \leq f(d)\} \\
=\ & \bigcup \{\mathcal{U}_d : y \leq f(d)\} \\
&
\end{align*}
is open in $I_{T}(X)$. For any $f: X \to Y$ and any $g: Y \to Z$, $T(g \circ f) (D) = \da g(f(D)) = \da g(\da f(D)) = T(g) (T(f) (D))$. Therefore, $T$ is a endofunctor. $\Box$

\begin{lem} \rm \label{product topological ideal}
${\rm I}_{T}(X \otimes Y) = {\rm I}_{T}(X) \otimes {\rm I}_{T}(Y)$
\end{lem}
\noindent{\bf Proof.} Since ${\rm I}_{T}(X)$ is algebraic, then it is core compact and ${\rm I}_{T}(X) \otimes {\rm I}_{T}(Y) = {\rm I}_{T}(X) \times {\rm I}_{T}(Y)$ by Theorem 
\ref{core compact}, where ${\rm I}_{T}(X) \times {\rm I}_{T}(Y)$ is the topological product of  ${\rm I}_{T}(X)$ and ${\rm I}_{T}(Y)$. It is easy to check that given any topological ideal $D \in {\rm I}_{T}(X \otimes Y)$, $D = \pi_1 D \times \pi_2 D$ with the property that $\pi_1 D, \pi_2 D$ are topological ideals in $X,Y$ respectively, and vice versa, where $\pi_1,\pi_2$ are the projections on $X,Y$ respectively. Then we need only to prove that their topologies are equal. Since $\mathcal{U}_x$ and $\mathcal{U}_y$ are bases in $X$ and $Y$ respectively, and  $\mathcal{U}_{(x,y)} = \{D \in {\rm I}_{T}(X \otimes Y) : (x,y) \in D\} = \mathcal{U}_x \times \mathcal{U}_y$ for any $x \in X, y\in Y$, the bases of topology of ${\rm I}_{T}(X) \times {\rm I}_{T}(Y)$ and topology $\Omega({\rm I}_T(X \otimes Y))$ are equal. Therefore, ${\rm I}_{T}(X \otimes Y) = {\rm I}_{T}(X) \otimes {\rm I}_{T}(Y)$. $\Box$

\begin{defn} \rm 
We define $\overline{T}$ as follows.

$$\forall \langle X,\{f_i\}_I \rangle \in {\rm Ob}({\bf DTop}(\Sigma,\mathcal{E})), \overline{T}(\langle X,\{f_i\}_I\rangle ) = \langle {\rm I}_T(X),\{\overline{f_i}\}_I\rangle,$$ where
 $$\overline{f_i}: {\rm I}_{T}(X)^n \to {\rm I}_{T}(X) \text{ with } \overline{f_i}(D_1,\dots,D_n) = \da \{f_i(d_1,\dots,d_n): d_j \in D_j\}.$$ 
And
 $$\forall f \in {\rm Mor}({\bf DTop}(\Sigma,\mathcal{E})), \overline{T}(f) = T(f). $$
\end{defn}

\begin{lem} \rm 
$\overline{T}$ is an endofunctor on ${\bf DTop}(\Sigma,\mathcal{E})$ such that $TU = U\overline{T}$.
\end{lem}
\noindent{\bf Proof.} 
{\bf Step 1.} $\overline{f_i}$ is well defined and continuous. $\langle T(X),\{\overline{f_i}\}_I \rangle$ satisfies $\mathcal{E}$.

Suppose that $D = D_1 \times \dots  \times D_{n_i} \in {\rm I}_{T}(X)^{n_i}$ and $D \to \ovr{x} = (x_1,x_2,\dots,x_{n_i}) = \bigvee D \in X^{n_i}$ .
Since $f_i : X^{n_i} \to X$ is continuous, we have $\{f_i(d_1,\dots,d_{n_i})\}_{d_j \in D_j} \to f_i(x_1,\dots,x_n)$ and then $\bigvee \{f_i(d_1,\dots,d_{n_i})\}_{d_j \in D_j} = f_i(x_1,\dots,x_{n_i}) $. Thus, $\overline{f_i}(D) \in {\rm I}_T(X)$. For continuity, consider $\overline{f_i}^{-1}(\big\uparrow_{{\rm I}_T(X)} (\da x))$. We have
\begin{align*} 
   & \overline{f_i}^{-1}(\big\uparrow_{{\rm I}_T(X)} (\da x)) \\
=\ & \{D \in  {\rm I}_{T}(X)^{n_i} : x \in \overline{f_i}(D) \}  \\ 
=\ & \{D \in  {\rm I}_{T}(X)^{n_i} : \exists (d_1,\dots,d_{n_i}) \in D.\ x \leq f_i(d_1,\dots,d_{n_i}) \} \\
=\ & \bigcup \{\mathcal{U}_{d_1} \times \cdots \times \mathcal{U}_{d_{n_i}}:\ x \leq f_i(d_1,\dots,d_{n_i})\}
&
\end{align*}
is open in ${\rm I}_T(X)^n$. Since $f_i \leq f_j$ implies $\overline{f_i} \leq \overline{f_j}$, it is easy to verify that inequalities in $\mathcal{E}$ hold for $\langle I_T(X),\{\overline{f_i}\}_I\rangle$ as well.
\vskip 3mm
\noindent{\bf Step 2.} $\overline{T}(f)$ is a homomorphism. Assume that $f:\langle X,\{f_i\}_I\rangle \to \langle Y,\{g_i\}_I\rangle$ is a continuous homomorphism. We need only to check the commutativity of $\overline{T}(f)$ and operations $\overline{f_i}$.
\begin{align*}
& \overline{T}(f) (\overline{f_i}(D_1,\dots,D_{n_i})) \\
=\ & \da f ( \overline{f_i}(D_1,\dots,D_{n_i}))\\
=\ & \da f( \da \{f_i(d_1,\dots,d_{n_i}): d_j \in D_j \}) \\
=\ & \da f( \{f_i(d_1,\dots,d_{n_i}): d_j \in D_j \})  \\
=\ &\da \{ g_i(f(d_1),\dots,f(d_{n_i})): d_j \in D_j \} \\
=\ & \da \{g_i(e_1,\dots,e_{n_i}) : e_j \in \da f(D_j)\} \\
=\ & \overline{g_i}(T(f)(D_1),\dots,T(f)(D_{n_i}) ). \\
&
\end{align*}
{\bf Step 3.} $TU= U\overline{T}$. Obviously. $\Box $

\begin{defn} \rm 
If a continuous map $f: X \to Y$ between continuous spaces satisfies $\da (f(\dda x)) \subseteq \dda f(x)$, i.e., $x \ll y$ implies $f(x) \ll f(y)$, we call $f$ way-below preserving.
\end{defn}

We denote $F$ the left adjoint of forgetful functor $U: {\bf DTop}(\Sigma,\mathcal{E}) \longrightarrow {\bf DTop}$. 

\begin{prop}\rm  \label{con way}
Let $X$ be a continuous space. Then $UF(X)$ is a continuous space and the unit $\eta: X \to UF(X) $ as well as the $\Sigma$-operations of $F(X)$ are way-below preserving. 
\end{prop}

$$\xymatrix{
  X \ar[r]^{\eta_X} \ar[d]_{\dda} &  UF(X) \ar[d]^{Uk_X}  & F(X) \ar[d]_{k_X}   \\
  T(X) \ar[r]_{T\eta_X} \ar[d]_{\sup_X}&      TUF(X)=U\overline{T}F(X) \ar[d]_{\sup_{UF(X)}} & \overline{T}F(X) \ar[d]_{\sup_{UF(X)}}  \\
  X \ar[r]_{\eta_X} & UF(X)  &  F(X) } $$

\noindent{\bf Proof.}
Since $\eta_X:X \to UF(X)$ has the universal property, there exists a unique continuous homomorphism $k_X$ such that the upper rectangle in the left of the diagram commutes. 

By the continuity of $\eta_X$, 
given any $D \in T(X) = I_{T}(X)$, $\eta_X \circ \sup_X (D) = \eta_X(\sup_X D) = \sup_{UF(X)}(\eta_X(D))$. Since $\sup_{UF(X)} \circ T\eta_X (D) = \sup_{UF(X)} (\da \eta_X(D) ) = \sup_{UF(X)}(\eta_X(D))$, the lower rectangle of the diagram commutes.

 Since $X$ is continuous, $\sup_X \circ \dda_X = 1_X$ and $\sup_{UF(X)}$ is a homomorphism from $\overline{T}F(X) $ to $F(X)$,  we have $\sup_{UF(X)} \circ Uk_X = 1_{UF(X)}$ by the universal property of $\eta_X$. Similarly, since $\dda_X \circ \sup_X \leq 1_{T(X)}$,  we have $Uk_X \circ  \sup_{UF(X)} \circ T\eta_X = T\eta_X \circ \dda_X \circ \sup_X   \leq T\eta_X $. Thus,  $Uk_X \circ  \sup_{UF(X)} \leq  1_{U\overline{T}F(X)}$. Then, $k_X = \dda_{UF(X)}$  and $F(X)$ is a continuous space by Theorem \ref{retract}. 

 We have $\dda_{(UF(X))}  \circ \eta_X (y) = \dda_{UF(X)}(\eta_X (y)) = T\eta_X \circ \dda_X (y) = \big\downarrow_{UF(X)} \eta_X(\dda y)$. Assume that $x \ll y$, then $\eta_X(x) \in \dda_{UF(X)} (\eta_X(y))$, i.e., $\eta_X(x) \ll_{UF(X)} \eta_X(y)$. Thus, 
 $\eta_X$ is way-below preserving.

By Lemma \ref{prod way}, the way-below map on a finite product of continuous spaces is the product of the individual way-below maps. Since the way-below map for $UF(X)$ lifts to a homomorphism from $F(X)$ to $\overline{T}F(X)$, it commutes with all $\Sigma$-operations $\{f_i\}_I$ of $F(X)$. We have
\begin{align*}
& \dda_{F(X)}(f_i(x_1,\dots,x_{n_i}))     \\
=\ & \overline{f_i}(\dda_{F(X)}(x_1),\dots,\dda_{F(X)}(x_{n_i}))  \\
=\ & \da {f_i}(\dda_{F(X)^{n_i}} (x_1,x_2,\dots,x_{n_i})).
\end{align*}
Therefore, these operations of $F(X)$ are way-below preserving. $\Box$

\begin{prop} \rm 
For every algebraic space $X$, the space $UF(X)$ is algebraic.
\end{prop}
\noindent{\bf Proof.} Let $S = \eta_X(X)$ be the subspace of $UF(X)$. Denote $\overline{S}$ the subalgebra of $UF(X)$ generated by $S$. Consider $\mathcal{D}\overline{S}$ together with the operations inherent from $F(X)$. Let $f$ be the same set map of $\eta_X $ from $X$ to $\mathcal{D}\overline{S}$. Then $f$ is continuous. By the universal property of $\eta_X$, there exists a $\overline{f}:F(X) \to \mathcal{D}\overline{S}$  such that $U\overline{f} \circ \eta_X = f$. Consider the following diagram.

$$\xymatrix{
  X \ar[r]^{\eta_X} \ar[rd]_{f} &  UF(X) \ar[d]^{U\overline{f}} \ar[rd]^{1_{UF(X)}}  \\
    & \mathcal{D}\overline{S} \ar[r]_{em} & UF(X)                   }$$

By the universal property of $\eta_X$, 
$U\overline{f} \circ \eta_X = f$. Considering $em \circ f :X \to UF(X)$, by the universal property of $\eta_X$, we have $em \circ U\overline{f} = 1_{UF(X)}$. Since $em$ is an injection, we have $em$ and $U\overline{f}$ are bijections and $\mathcal{D}\overline{S} = UF(X)$. Assume that $x = f_i(x_1,\dots,x_{n_i})$ and $x_j = \eta(a_j)\in S, a_j \in X$ for each $1 \leq j \leq n_i$ and $f_i$ an operation of $F(X)$. Since $X$ is algebraic, for each $a_j$, there is a directed subset $D_j \subseteq K(X)$ such that $D_j \leq a_j$ and $D_j \to a_j$. Since $\eta_X$ is continuous and way-below preserving, we have $\eta(D_j) \to x_j$, $\eta(D_j) \leq x_j$ and $\eta(D_j) \subseteq K(UF(X))$. Since $f_i$ is way-below preserving, we have that $f_i(\eta(D_1) \times \cdots \times \eta( D_n) ) \subseteq K(UF(X)) $ and it converges to $f_i(x_1,\dots,x_n)$. 

Denote $S_0 = S$ and $S_{n+1} = S_n \cup \{f_i(x_1,\dots,x_{n_i}) : i \in I\  \&\ x_j \in S_n \text{ for } 1 \leq j \leq n_i\}$ for $n \in \mathbb{N}$. Then $UF(X) = \overline{S} = \bigcup_{n \in \mathbb{N}} S_n$.
The same, assume that $y_j \in S_n$ and there exists a directed set $D_j \subseteq K(UF(X))$ such that $D_j \to x_j$ and $D_j \leq x_j$ for $1 \leq j \leq n_i$. Then $f_i(D_1\times \dots \times D_n)$ is a directed subset of compact elements of $UF(X)$ that is lower than and converges to $ f_i(x_1,\dots, x_{n_i}) \in S_{n+1}$. 

Thus, by structure induction, we know that for each element $x$ of $\overline{S}$, there is a directed subset of compact elements of $UF(X)$ that is lower than and converges to $x$. Thus, $UF(X)$ is algebraic. $UF(X) = \mathcal{D}\overline{S} = \overline{S}$. $\Box$

\end{document}